\documentclass[12pt]{article}
\usepackage{graphicx}
\usepackage{amsmath,amsthm,amssymb,enumerate}%, esint}
\usepackage{euscript,mathrsfs}
\usepackage{color}
\usepackage{dsfont}
\usepackage[left=2cm,right=2cm,top=3.5cm,bottom=3.5cm]{geometry}
\usepackage{color}
\usepackage[framemethod=tikz]{mdframed}
\allowdisplaybreaks

\newcommand{\tc}{\color{black}}

\usepackage{soul}
\usepackage{esint}
\catcode`\@=11 \@addtoreset{equation}{section}

\catcode`\@=12

\allowdisplaybreaks

\newtheorem{Theorem}{Theorem}[section]
\newtheorem{Proposition}[Theorem]{Proposition}
\newtheorem{Lemma}[Theorem]{Lemma}
\newtheorem{Corollary}[Theorem]{Corollary}

\theoremstyle{definition}
\newtheorem{Definition}[Theorem]{Definition}

\newtheorem{Remark}[Theorem]{Remark}

\newcommand{\bTheorem}[1]{
\begin{Theorem} \label{T#1} }
\newcommand{\eT}{\end{Theorem}}

\newcommand{\bProposition}[1]{
\begin{Proposition} \label{P#1}}
\newcommand{\eP}{\end{Proposition}}

\newcommand{\bLemma}[1]{
\begin{Lemma} \label{L#1} }
\newcommand{\eL}{\end{Lemma}}

\newcommand{\bCorollary}[1]{
\begin{Corollary} \label{C#1} }
\newcommand{\eC}{\end{Corollary}}

\newcommand{\bRemark}[1]{
\begin{Remark} \label{R#1} }
\newcommand{\eR}{\end{Remark}}

\newcommand{\bDefinition}[1]{
\begin{Definition} \label{D#1} }
\newcommand{\eD}{\end{Definition}}

\newcommand{\Del}{\Delta_x}

\newcommand{\Ds}{\mathbb{D}_x}

\newcommand{\MTC}{\mathcal{T}}

\newcommand{\tvm}{\widetilde{\vc{m}}}
\newcommand{\tS}{\widetilde{S}}
\newcommand{\bfphi}{\boldsymbol{\varphi}}

\newcommand{\bFormula}[1]{
\begin{equation} \label{#1}}
\newcommand{\eF}{\end{equation}}

\newcommand{\Ov}[1]{\overline{#1}}

\newcommand{\aleq}{\lesssim} %changed by Florian

\newcommand{\ageq}{\stackrel{>}{\sim}}

\newcommand{\vr}{\varrho}
\newcommand{\vre}{\vr_\ep}

\newcommand{\vte}{\vt_\ep}
\newcommand{\vue}{\vu_\ep}
\newcommand{\tvr}{\tilde \vr}
\newcommand{\tvu}{{\tilde \vu}}
\newcommand{\tvt}{\tilde \vt}

\newcommand{\vt}{\vartheta}
\newcommand{\vu}{\vc{u}}
\newcommand{\vm}{\vc{m}}

\newcommand{\vc}[1]{{\bf #1}}

\newcommand{\Div}{{\rm div}_x}
\newcommand{\Grad}{\nabla_x}

\newcommand{\dx}{\,{\rm d} {x}}

\newcommand{\dt}{\,{\rm d} t }

\newcommand{\vU}{\vc{U}}

\newcommand{\intO}[1]{\int_{\Omega} #1 \ \dx}
\newcommand{\avintO}[1]{\fint_{\Omega} #1 \dx}

\newcommand{\D}{{\rm d}}

\newcommand{\ep}{\varepsilon}

\newcommand{\R}{\mathbb{R}}
\newcommand{\vtB}{\vt_B}

\newcommand{\br}{ \nonumber \\ }

\newcommand{\cred}{\color{black}} %showing/hiding comments

\def\softd{{\leavevmode\setbox1=\hbox{d}%
          \hbox to 1.05\wd1{d\kern-0.4ex{\char039}\hss}}}
\definecolor{Cgrey}{rgb}{0.85,0.85,0.85}
\definecolor{Cblue}{rgb}{0.50,0.85,0.85}
\definecolor{Cred}{rgb}{1,0,0}
\definecolor{fancy}{rgb}{0.10,0.85,0.10}

\newcommand\Cbox[2]{%
    \newbox\contentbox%
    \newbox\bkgdbox%
    \setbox\contentbox\hbox to \hsize{%
        \vtop{
            \kern\columnsep
            \hbox to \hsize{%
                \kern\columnsep%
                \advance\hsize by -2\columnsep%
                \setlength{\textwidth}{\hsize}%
                \vbox{
                    \parskip=\baselineskip
                    \parindent=0bp
                    #2
                }%
                \kern\columnsep%
            }%
            \kern\columnsep%
        }%
    }%
    \setbox\bkgdbox\vbox{
        \color{#1}
        \hrule width  \wd\contentbox %
               height \ht\contentbox %
               depth  \dp\contentbox
        \color{black}
    }%
    \wd\bkgdbox=0bp%
    \vbox{\hbox to \hsize{\box\bkgdbox\box\contentbox}}%
    \vskip\baselineskip%
}

\mdfdefinestyle{MyFrame}{%
	linecolor=black,
	outerlinewidth=1pt,
	roundcorner=5pt,
	innertopmargin=\baselineskip,
	innerbottommargin=\baselineskip,
	innerrightmargin=10pt,
	innerleftmargin=10pt,
	backgroundcolor=white!20!white}

\date{}

\begin{document}

%%%%%%%%%%%%%%%%%%%%%%%%%%%%%%%%

\title{Rigorous derivation of the Oberbeck-Boussinesq approximation revealing unexpected term}

\author{Peter Bella$^{1,}$\thanks{P.B. and F.O. were partially supported by the German {\tc Research} Foundation DFG in context of the Emmy Noether Junior Research Group BE 5922/1-1.}
\and Eduard Feireisl$^{2,}$\thanks{The work of E.F. was partially supported by the
Czech Sciences Foundation (GA\v CR), Grant Agreement
21--02411S. The Institute of Mathematics of the Academy of Sciences of
the Czech Republic is supported by RVO:67985840. } \and Florian Oschmann$^{1,\ast}$
}

\date{}

\maketitle

\medskip

\centerline{$^1$ TU Dortmund, Fakult\"at f\"ur Mathematik}

\centerline{Vogelpothsweg 87, 44227 Dortmund, Germany}

\medskip

\centerline{$^2$ Institute of Mathematics of the Academy of Sciences of the Czech Republic}

\centerline{\v Zitn\' a 25, CZ-115 67 Praha 1, Czech Republic}

\date{}

\maketitle

\begin{abstract}
	
	We consider a general compressible viscous and heat conducting fluid confined between two parallel plates 
	and heated from the bottom. The time evolution of the fluid is described by the Navier--Stokes--Fourier system considered in the regime of low Mach and Froude numbers suitably interrelated. Surprisingly and differently to the case of Neumann boundary conditions for the temperature, the asymptotic limit is identified  as the Oberbeck--Boussinesq system supplemented with non--local boundary conditions for the temperature deviation. 
	
\end{abstract}

%\bigskip


% \tableofcontents

\section{Introduction}
\label{i}

The Rayleigh--B\' enard problem concerns the motion of a compressible fluid confined between two parallel plates 
heated from the bottom. For the sake of simplicity, we suppose the motion is space--periodic with respect to the horizontal variable. Accordingly, the underlying spatial domain may be identified with the flat torus 
in the horizontal plane:
\begin{equation*}
\Omega = \mathbb{T}^{d - 1}  \times (0,1),\ \mathbb{T}^{d-1} = \left( [0,1] \Big|_{\{ 0,1 \} } \right)^{d - 1},\ d=2,3. 
\end{equation*}
The state of the fluid at a given time $t$ and a spatial position $x \in \Omega$ is characterized by the 
mass density $\vr = \vr(t,x)$, the absolute temperature $\vt = \vt(t,x)$, and the velocity $\vu = \vu(t,x)$. 
We consider the motion in the asymptotically incompressible and stratified regime. Accordingly, 
the time evolution of the flow is governed by the scaled \emph{Navier--Stokes--Fourier (NSF) system}:

\begin{mdframed}[style=MyFrame]
 
\begin{align} 
	\partial_t \vr + \Div (\vr \vu) &= 0, \label{i1}\\
	\partial_t (\vr \vu) + \Div (\vr \vu \otimes \vu) + \frac{1}{\ep^2} \Grad p(\vr, \vt) &= \Div \mathbb{S}(\vt, \Grad \vu) + \frac{1}{\ep} \vr \Grad G, \label{i2} \\ 
	\partial_t (\vr s(\vr, \vt)) + \Div (\vr s (\vr, \vt) \vu) + \Div \left( \frac{ \vc{q} (\vt, \Grad \vt) }{\vt} \right) &= 
	\frac{1}{\vt} \left( \ep^2 \mathbb{S} : \Grad \vu - \frac{\vc{q} (\vt, \Grad \vt) \cdot \Grad \vt }{\vt} \right),
	\label{i3}	
\end{align}
supplemented with the boundary conditions 
\begin{align} 
	\vu|_{\partial \Omega} &= 0, \label{i4} \\
	\vt|_{\partial \Omega} &= \vtB.
	\label{i5}
	\end{align}

\end{mdframed}
The viscous stress is given by Newton's rheological law 
\begin{equation} \label{i6}
	\mathbb{S}(\vt, \Grad \vu) = \mu(\vt) \left( \Grad \vu + \Grad^t \vu - \frac{2}{d} \Div \vu \mathbb{I} \right) + \eta(\vt) \Div \vu \mathbb{I}, 
	\end{equation}
while the internal energy flux satisfies Fourier's law
\begin{equation} \label{i7} 
	\vc{q}(\vt, \Grad \vt) = - \kappa (\vt) \Grad \vt.
	\end{equation}

Here, the more conventional internal energy (heat) equation
\begin{equation} \label{i8}
	\partial_t (\vr e(\vr, \vt)) + \Div (\vr e(\vr, \vt) \vu) + \Div \vc{q} (\vt, \Grad \vt) = 
	\ep^2 \mathbb{S}(\vt, \Grad \vu): \Grad \vu - p (\vr, \vt) \Div \vu
	\end{equation}
is replaced by the entropy balance equation \eqref{i3}. {\tc We point out that the equations 
\eqref{i3}, \eqref{i8} are \emph{equivalent} 
in the framework of classical solutions as long as the thermodynamic functions $p(\vr, \vt)$, $e(\vr, \vt)$, 
$s(\vr, \vt)$ are interrelated through Gibbs' equation
\begin{equation} \label{i9} 
	\vt D s = D e + p D \left( \frac{1}{\vr} \right).
	\end{equation}
The entropy balance equation, however, is more convenient for the weak formulation introduced in Section \ref{L}. In particular, 
the weak solutions based on the entropy formulation enjoy several useful properties, notably the so--called weak strong uniqueness principle.
}

Besides Gibbs' equation, we impose the hypothesis of \emph{thermodynamic stability} written in the form 
\begin{equation} \label{HTS}
	\frac{\partial p(\vr, \vt) }{\partial \vr } > 0,\ 
	\frac{\partial e(\vr, \vt) }{\partial \vt } > 0 \ \mbox{for all}\ \vr, \vt > 0.
\end{equation}

The scaling in \eqref{i1}--\eqref{i3} represented by a small parameter $\ep > 0$ corresponds to the 
Mach number ${\rm Ma} = \ep$ (the ratio between a ``characteristic'' fluid's velocity and the local speed of sound) and the Froude number ${\rm Fr} = \ep^{\frac{1}{2}}$ (the ratio between inertial and gravitational forces of the fluid), see e.g.~the survey paper by Klein et al.~\cite{KBSMRMHS}.
Our goal is to identify the asymptotic limit of solutions for $\ep \to 0$.

\subsection{Asymptotic limit}

The zero--th order terms in the asymptotic limit are determined by the stationary problem 
\begin{equation} \label{i10}
	\Grad p(\vr, \vt ) = \ep \vr \Grad G.
\end{equation} 
Accordingly, a solution of the Navier--Stokes--Fourier system can be written in the form 
\begin{equation*}
\vre = \Ov{\vr} + \ep \mathfrak{R}_\ep,\ \vte = \Ov{\vt} + \ep \mathfrak{T}_\ep,
\end{equation*}
where $\Ov{\vr}$, $\Ov{\vt}$ are positive constants and $\mathfrak{R}_\ep \to \mathfrak{R}$, 
$\mathfrak{T}_\ep \to \mathfrak{T}$ satisfy in the asymptotic limit $\ep \to 0$ the so called \emph{Boussinesq relation}
\begin{equation} \label{i11}
\frac{\partial p(\Ov{\vr}, \Ov{\vt} ) }{\partial \vr} \Grad \mathfrak{R} + 
\frac{\partial p(\Ov{\vr}, \Ov{\vt} ) }{\partial \vt} \Grad \mathfrak{T}	 = \Ov{\vr} \Grad G.
	\end{equation}
Without loss of generality, we suppose
\begin{equation*}
\intO{ G } = 0,\   \Ov{\vr} = \avintO{ \vre } \equiv \frac{1}{|\Omega|} \intO{\vre } 
\ \Rightarrow\ \intO{ \mathfrak{R}_\ep } = 0,
\end{equation*}
meaning, in particular, the total mass of the fluid is constant independent of $\ep$.

Anticipating convergence of the temperature deviations 
\begin{equation*}
\mathfrak{T}_\ep = \frac{\vte - \Ov{\vt}}{\ep} \to \mathfrak{T} 
\end{equation*}
we impose the boundary condition 
\begin{equation} \label{i11a}
\vte|_{\partial \Omega} = \vt_B = \Ov{\vt} + \ep \Theta_B,\ \Theta_B = \Theta_B (x),
\end{equation}
where the perturbation $\Theta_B$ is not necessarily positive.

\subsection{Limit system}

Formally, the limit problem is expected to be the \emph{Oberbeck--Boussinesq (OB) system}:
\begin{align} 
	\Div \vU &= 0, \br	
	\Ov{\vr} \Big( \partial_t \vU + \vU \cdot \Grad \vU \Big) + \Grad \Pi &= \Div \mathbb{S}(\Ov{\vt}, \Grad \vU) +  \mathfrak{R} \Grad G, \br
	\Ov{\vr} c_p(\Ov{\vr}, \Ov{\vt} ) \left( \partial_t \mathfrak{T} + \vU \cdot \Grad \mathfrak{T} \right)	- 
	\Ov{\vr} \ \Ov{\vt} \alpha(\Ov{\vr}, \Ov{\vt} ) \vU \cdot \Grad G
	&= \kappa(\Ov{\vt}) \Del \mathfrak{T},	
	\label{ObBs}
\end{align}
where 
\begin{align}\label{Coeff} 
	\alpha(\Ov{\vr}, \Ov{\vt} ) &\equiv \frac{1}{\Ov{\vr}}  \frac{\partial p(\Ov{\vr}, \Ov{\vt} ) }{\partial \vt} \left( \frac{\partial p(\Ov{\vr}, \Ov{\vt} ) }{\partial \vr} \right)^{-1}, \br 
	c_p (\Ov{\vr}, \Ov{\vt} ) &\equiv \frac{\partial e(\Ov{\vr}, \Ov{\vt} ) }{\partial \vt}	+ \Ov{\vr}^{-1} \Ov{\vt}
	\alpha(\Ov{\vr}, \Ov{\vt} ) \frac{\partial p(\Ov{\vr}, \Ov{\vt} ) }{\partial \vt}
	%\nonumber
\end{align}
stand for the coefficient of thermal expansion,
and the specific heat at constant pressure, respectively. The quantities $\mathfrak{R}$, $\mathfrak{T}$ 
satisfy the Boussinesq relation \eqref{i11}.

We refer to the survey by Zeytounian \cite{ZEY1} and the list of references therein for a formal derivation of the Oberbeck--Boussinesq system. A rigorous proof of convergence of a family of (weak) solutions $(\vre, \vte, \vue)_{\ep > 0}$ to 
the Navier--Stokes--Fourier system, meaning 
\begin{equation*}
\vue \to \vU , \ \frac{\vte - \Ov{\vt}}{\ep} \to \mathfrak{T} \ \mbox{in some sense,}
\end{equation*}
was obtained in \cite{FN5} (see also \cite[Chapter~5]{FeNo6A}) for the conservative boundary conditions 
\begin{align} 
	\vue \cdot \vc{n}|_{\partial \Omega} &= 0,\ 
	(\mathbb{S}(\vte, \Grad \vue) \cdot \vc{n}) \times \vc{n}|_{\partial \Omega} = 0, \br 
	\Grad \vte \cdot \vc{n}|_{\partial \Omega} &= 0.
	\label{i16}	
\end{align}
Note that the conservative boundary conditions \eqref{i16} imply that the average
\begin{equation*}
\avintO{ \mathfrak{T} (t, \cdot) }
\end{equation*}
is a constant of motion. In particular, the Boussinesq relation \eqref{i11} can be written as
\begin{equation} \label{i16AA}
\frac{\partial p(\Ov{\vr}, \Ov{\vt} ) }{\partial \vr} \mathfrak{R} + 
\frac{\partial p(\Ov{\vr}, \Ov{\vt} ) }{\partial \vt} \mathfrak{T}	 = \Ov{\vr} G
\end{equation}
recalling the convention $\intO{G } = 0$. The form \eqref{i16AA} of the Boussinesq relation is rather misleadingly interpreted in the sense that the Oberbeck--Boussinesq approximation can be 
derived from the Navier--Stokes--Fourier system on condition that the density is supposed to be an affine function of the temperature, cf.~the discussion in Maruyama \cite{maruyama2019rational}.

{\tc
The present paper can be seen as the first attempt to address the same singular limit problem under the inhomogeneous Dirichlet boundary conditions 
imposed on the temperature, physically relevant for the well--known Rayleigh--B\' enard convection problem. The result is stated in the framework of the 
general theory of weak solutions to the primitive Navier--Stokes--Fourier system with inhomogeneous boundary conditions developed recently in \cite{ChauFei} and the 
monograph \cite{FeiNovOpen}. }   

Rather surprisingly, we show that the limit system is in fact \emph{different} if the Dirichlet boundary conditions \eqref{i4}, \eqref{i5} are imposed. In accordance with \eqref{i4}, \eqref{i11a}, we consider the 
scaled Navier--Stokes--Fourier system with the 
boundary conditions
\begin{equation} \label{scb}
\vue|_{\partial \Omega}= 0,\ \vte|_{\partial \Omega} = \vtB = \Ov{\vt} + \ep \Theta_B.	
\end{equation}
Imposing \eqref{scb} in place of \eqref{i16} drastically changes the behaviour of the fluid 
as the resulting system is no longer energetically closed and the dynamics is driven by the temperature gradient. Indeed the total mass conservation discussed by several authors 
(Barletta et al. \cite{barletta2022Bous,barletta2022use}, Maruyama \cite{maruyama2019rational}) must be enforced through the 
condition
\begin{equation*}
\intO{ \mathfrak{R}_\ep } = 0 \ \Rightarrow \ \intO{ \mathfrak{R} } = 0
\end{equation*}
in \eqref{i11}. Accordingly, the conventional Boussinesq relation \eqref{i16AA} must be replaced by
\begin{equation} \label{OBchi}
\frac{\partial p(\Ov{\vr}, \Ov{\vt} ) }{\partial \vr} \mathfrak{R} + 
\frac{\partial p(\Ov{\vr}, \Ov{\vt} ) }{\partial \vt} \mathfrak{T}	 = \Ov{\vr} G + \chi(t),
\end{equation}
where, in accordance with our convention $\intO{G} = 0$,
\begin{equation} \label{relB}
\chi = \frac{\partial p(\Ov{\vr}, \Ov{\vt} ) }{\partial \vt}  \avintO{ \mathfrak{T} }.
\end{equation}

Relation \eqref{relB} indicates that the limit problem could involve ``non--local'' terms. 
Indeed we show that a family of (weak) solutions to the NSF system $(\vre, \vte, \vue)_{\ep > 0}$ admits a limit

\begin{equation*}
	\vue \to \vU , \ \frac{\vte - \Ov{\vt}}{\ep} \to \mathfrak{T},\ \frac{\vre - \Ov{\vr}}{\ep} \to \mathfrak{R}
\end{equation*}
solving a \emph{modified Oberbeck--Bousinesq system}
\begin{align} 
	\Div \vU &= 0, \br	
	\Ov{\vr} \Big( \partial_t \vU + \vU \cdot \Grad \vU \Big) + \Grad \Pi &= \Div \mathbb{S}(\Ov{\vt}, \Grad \vU) +  \mathfrak{R} \Grad G, \br
	\Ov{\vr} c_p(\Ov{\vr}, \Ov{\vt} ) \left( \partial_t \mathfrak{T} + \vU \cdot \Grad \mathfrak{T} \right)	- 
	\Ov{\vr} \ \Ov{\vt} \alpha(\Ov{\vr}, \Ov{\vt} ) \vU \cdot \Grad G
	&= \kappa(\Ov{\vt}) \Del \mathfrak{T} + \Ov{\vt} \alpha (\Ov{\vr}, \Ov{\vt}) \frac{\partial p (\Ov{\vr}, \Ov{\vt})}{\partial \vt}  
	\partial_t \avintO{ \mathfrak{T} },	
	\label{OBsI}
\end{align}
together with the Boussinesq relation
\begin{equation} \label{OB1I}
	\frac{\partial p(\Ov{\vr}, \Ov{\vt} ) }{\partial \vr} \Grad \mathfrak{R} + 
	\frac{\partial p(\Ov{\vr}, \Ov{\vt} ) }{\partial \vt} \Grad \mathfrak{T} = \Ov{\vr} \Grad G,\ 
	\intO{ \mathfrak{R} } = 0,
\end{equation} 
and the boundary conditions
\begin{equation} \label{bcI}
	\vU|_{\partial \Omega} = 0,\ \mathfrak{T}|_{\partial \Omega} = \Theta_B.
\end{equation}

To simplify even more, using relation~\eqref{OBchi} to replace $\mathfrak{R}$ with $\mathfrak{T}$ in the momentum equation while observing that two gradient terms hide into $\Pi$, the system~\eqref{OBsI}--\eqref{bcI} turns into 
\begin{equation}\label{P1}
\begin{aligned} 
	\Div \vU &= 0, \\
	\Ov{\vr} \Big( \partial_t \vU + \vU \cdot \Grad \vU \Big) + \Grad \Pi &= \Div \mathbb{S}(\Ov{\vt}, \Grad \vU) - \Ov{\vr} \alpha( \Ov{\vr},\Ov{\vt}) \mathfrak{T} \Grad G, \\
	\Ov{\vr} c_p(\Ov{\vr}, \Ov{\vt} ) \left( \partial_t \mathfrak{T} + \vU \cdot \Grad \mathfrak{T} \right)	- 
	\Ov{\vr} \ \Ov{\vt} \alpha(\Ov{\vr}, \Ov{\vt} ) \vU \cdot \Grad G
	&= \kappa(\Ov{\vt}) \Del \mathfrak{T} + \Ov{\vt} \alpha (\Ov{\vr}, \Ov{\vt}) \frac{\partial p (\Ov{\vr}, \Ov{\vt})}{\partial \vt}  
	\partial_t \avintO{ \mathfrak{T} }, \\
	\vU|_{\partial \Omega} &= 0,\ \mathfrak{T}|_{\partial \Omega} = \Theta_B. 
\end{aligned}
\end{equation}

Finally, performing a simple change of variables
\[
r = \mathfrak{R}, \ \Theta = \mathfrak{T} - \lambda(\Ov{\vr}, \Ov{\vt}) \avintO{ \mathfrak{T} }, \ \mbox{where}\ \lambda(\Ov{\vr}, \Ov{\vt}) = \frac{\Ov{\vt} \alpha (\Ov{\vr}, \Ov{\vt} ) }
{\Ov{\vr} c_p(\Ov{\vr}, \Ov{\vt}) } \frac{\partial p(\Ov{\vr}, \Ov{\vt} )}{\partial \vt} \in (0,1),
\]
we may rewrite system \eqref{OBsI}--\eqref{bcI} in the form of \emph{conventional Oberbeck--Boussinesq system}
\begin{mdframed}[style=MyFrame]
	\begin{align} 
		\Div \vU &= 0, \br	
		\Ov{\vr} \Big( \partial_t \vU + \vU \cdot \Grad \vU \Big) + \Grad \Pi &= \Div \mathbb{S}(\Ov{\vt}, \Grad \vU) +  r \Grad G, \br
		\Ov{\vr} c_p(\Ov{\vr}, \Ov{\vt} ) \Big( \partial_t \Theta + \vU \cdot \Grad \Theta \Big)	- 
		\Ov{\vr} \ \Ov{\vt} \alpha(\Ov{\vr}, \Ov{\vt} ) \vU \cdot \Grad G
		&= \kappa(\Ov{\vt}) \Del \Theta, \br
		\frac{\partial p(\Ov{\vr}, \Ov{\vt} ) }{\partial \vr} \Grad r + 
		\frac{\partial p(\Ov{\vr}, \Ov{\vt} ) }{\partial \vt} \Grad \Theta	 &= \Ov{\vr} \Grad G	
		\label{ObBs1}
	\end{align}
	with the no--slip boundary condition for the velocity
	\begin{equation} \label{i18bis}
		\vU|_{\partial \Omega} = 0,
	\end{equation}
	and a \emph{non--local} boundary condition for the temperature deviation
	\begin{equation} \label{i14B}
		\Theta|_{\partial \Omega} = \Theta_B - \frac{\lambda(\Ov{\vr}, \Ov{\vt})}{1 - \lambda(\Ov{\vr}, \Ov{\vt})} \avintO{ \Theta }.	
	\end{equation}
\end{mdframed}
Here and hereafter, we always suppose that the coefficient of thermal expansion $\alpha (\Ov{\vr}, \Ov{\vt})$ at the reference state  $(\Ov{\vr}, \Ov{\vt})$ is strictly positive. In particular, 
it follows directly from \eqref{Coeff} and \eqref{HTS} that $\lambda\in (0,1)$. Parabolic equations with the non--local boundary terms similar to \eqref{i14B} have been recognized by Day \cite{Day1982} in the context of some models in thermoelasticity, and subsequently studied 
by a number of authors, see Chen and Liu \cite{ChenLiu}, Friedman \cite{Friedman86}, Gladkov and Nikitin \cite{GlaNik}, Pao \cite{Pao},  among others.

We point out that our results justify the Oberbeck--Boussinesq approximation as a rigorous description of the limit behaviour of solutions to the complete Navier--Stokes--Fourier system only in the regime where the temperature deviation from the constant state $\Ov{\vt}$ is small of order $\ep$. Indeed the relevance of the Oberbeck--Boussinesq approximation for the Rayleigh--B{\' e}nard convection with large deviation of the boundary temperature
has been questioned by several authors, see Bormann \cite{Borm}, Fr\" ohlich, Laure, and Peyret \cite{FrLaPe}, Nadolin \cite{Nadol}, among others. 

\subsection{The strategy of the convergence proof}

Our approach is based on the concept of \emph{weak solutions} to the NSF system with energetically 
open boundary conditions developed recently in \cite{ChauFei}, \cite{FeiNovOpen}. In particular, the relative energy inequality 
based on the ballistic energy balance is used to measure the distance between the solutions of the primitive and target systems. In contrast with \cite{FN5}, \cite[Chapter~5]{FeNo6A}, strong convergence is obtained 
with certain explicit estimates of the error depending on how close are the initial data to their limit values. 

The paper is organized as follows.
We start with the concept of weak solutions for the NSF system with Dirichlet boundary conditions 
introduced in \cite{ChauFei}. In particular, we recall the ballistic energy and the associated relative energy 
inequality in Section~\ref{L}. Then, in Section~\ref{OB}, we record the available results on solvability of the OB system in the framework of strong solutions. 
The main results on convergence to the target OB system are stated in Section~\ref{M}. The rest of the paper is then devoted to the proof of convergence.
In Section~\ref{B}, we derive the basic energy estimates that control the amplitude of the fluid velocity as well as 
the distance of the density and the temperature profiles from their limit values independent of the scaling 
parameter $\ep$. The proof of convergence to the OB system is completed in Section~\ref{WPL} by means of the 
relative energy inequality.

\section{Weak solutions to the primitive system}
\label{L}

Following \cite{ChauFei}, \cite{FeiNovOpen}, we introduce the concept of weak solutions to the NSF system. 

\begin{Definition}[{\bf Weak solution to the NSF system}] \label{DL1}
	We say that a trio $(\vr, \vt, \vu)$ is a weak solution of the NSF system \eqref{i1}--\eqref{i7}, 
	with the initial data
	\[
	\vr(0, \cdot) = \vr_0,\ (\vr \vu) (0, \cdot) = \vr_0 \vu_0,\ 
	(\vr s)(0, \cdot) = \vr_0 s(\vr_0, \vt_0),
	\]
if the following holds:	

\begin{itemize}
	
	\item The solution belongs to the \emph{regularity class}: 
	\begin{align}
		\vr &\in L^\infty(0,T; L^\gamma(\Omega)) \ \mbox{for some}\ \gamma > 1,\ \vr \geq 0 
		\ \mbox{a.a.~in}\ (0,T) \times \Omega, \br
		\vu &\in L^2(0,T; W^{1,2}_0 (\Omega; \mathbb{R}^d)), \br 
		\vt^{\beta/2} ,\ \log(\vt) &\in L^2(0,T; W^{1,2}(\Omega)) \ \mbox{for some}\ \beta \geq 2,\ 
		\vt > 0 \ \mbox{a.a.~in}\ (0,T) \times \Omega, \br
		(\vt - \vtB) &\in L^2(0,T; W^{1,2}_0 (\Omega)).
		\label{Lw6}
	\end{align}
	
	\item The \emph{equation of continuity} \eqref{i1} is satisfied in the sense of distributions, 
	more specifically,
	\begin{align} 
		\int_0^T \intO{ \Big[ \vr \partial_t \varphi + \vr \vu \cdot \Grad \varphi \Big] } \dt &=  - 
		\intO{ \vr(0) \varphi(0, \cdot) }
		\label{Lw4}
	\end{align}
	for any $\varphi \in C^1_c([0,T) \times \Ov{\Omega} )$.
	\item The \emph{momentum equation} \eqref{i2} is satisfied in the sense of distributions, 
	\begin{align}
		\int_0^T &\intO{ \left[ \vr \vu \cdot \partial_t \bfphi + \vr \vu \otimes \vu : \Grad \bfphi + 
			\frac{1}{\ep^2} p(\vr, \vt) \Div \bfphi \right] } \dt \br &= \int_0^T \intO{ \left[ \mathbb{S}(\vt, \Grad \vu) : \Grad \bfphi - \frac{1}{\ep} \vr \Grad G \cdot \bfphi \right] } \dt - 
		\intO{ \vr_0 \vu_0 \cdot \bfphi (0, \cdot) }
		\label{Lw5}
	\end{align}	
	for any $\bfphi \in C^1_c([0, T) \times \Omega; \mathbb{R}^d)$.
	
	\item The \emph{entropy balance} \eqref{i3} is replaced by the inequality
	\begin{align}
		- \int_0^T &\intO{ \left[ \vr s(\vr, \vt) \partial_t \varphi + \vr s (\vr ,\vt) \vu \cdot \Grad \varphi + \frac{\vc{q} (\vt, \Grad \vt )}{\vt} \cdot 
			\Grad \varphi \right] } \dt \br &\geq \int_0^T \intO{ \frac{\varphi}{\vt} \left[ \ep^2 \mathbb{S}(\vt, \Grad \vu) : \Grad \vu - 
			\frac{\vc{q} (\vt, \Grad \vt) \cdot \Grad \vt }{\vt} \right] } \dt + \intO{ \vr_0 s(\vr_0, \vt_0) 
			\varphi (0, \cdot) } 
		\label{Lw7} 
	\end{align}
	for any $\varphi \in C^1_c([0, T) \times \Omega)$, $\varphi \geq 0$.
	
	\item  The \emph{ballistic energy balance}
	\begin{align}  
		- \int_0^T &\partial_t \psi	\intO{ \left[ \ep^2 \frac{1}{2} \vr |\vu|^2 + \vr e(\vr, \vt) - \tvt \vr s(\vr, \vt) \right] } \dt  \br &+ \int_0^T \psi
		\intO{ \frac{\tvt}{\vt}	 \left[ \ep^2 \mathbb{S}(\vt, \Grad \vu): \Grad \vu - \frac{\vc{q}(\vt, \Grad \vt) \cdot \Grad \vt }{\vt} \right] } \dt  \br
		&\leq 
		\int_0^T \psi \intO{ \left[ \ep \vr \vu \cdot \Grad G - 
			\vr s(\vr, \vt) \partial_t \tvt - \vr s(\vr, \vt) \vu \cdot \Grad \tvt - \frac{\vc{q}(\vt, \Grad \vt)}{\vt} \cdot \Grad \tvt \right] }\br 
		&+ \psi(0) \intO{ \left[ \frac{1}{2} \ep^2 \vr_0 |\vu_0|^2 + \vr_0 e(\vr_0, \vt_0) - \tvt(0, \cdot) \vr_0 s(\vr_0, \vt_0) \right] }
		\label{Lw8}
	\end{align}
	holds for any $\psi \in C^1_c[0; T)$, $\psi \geq 0$, and any $\tvt \in C^1([0; T) \times \Ov{\Omega})$ satisfying
	\[
	\tvt > 0,\ \tvt|_{\partial \Omega} = \vtB.
	\]
\end{itemize}
 
	\end{Definition}

Although quite general, the weak solutions in the sense of Definition \ref{DL1} comply with the weak--strong uniqueness principle, meaning they coincide with the strong solution as long as the latter exists, see \cite[Chapter~4]{FeiNovOpen}.

\subsection{Relative energy inequality}

Following \cite{ChauFei}, \cite{FeiNovOpen}, we introduce the scaled \emph{relative energy} 
\begin{align}
	E_\ep &\left( \vr, \vt, \vu \Big| \tvr , \tvt, \tvu \right) \br &= \frac{1}{2}\vr |\vu - \tvu|^2 + 
	\frac{1}{\ep^2} \left[ \vr e - \tvt \Big(\vr s - \tvr s(\tvr, \tvt) \Big)- 
	\Big( e(\tvr, \tvt) - \tvt s(\tvr, \tvt) + \frac{p(\tvr, \tvt)}{\tvr} \Big)
	(\vr - \tvr) - \tvr e (\tvr, \tvt) \right] .
	\nonumber
\end{align}
The hypothesis of thermodynamic stability \eqref{HTS} can be equivalently rephrased as (strict) convexity 
of the total energy expressed with respect to the conservative entropy variables
\[
E_\ep \Big( \vr, S = \vr s(\vr, \vt), \vm = \vr \vu \Big) \equiv \frac{1}{2} \frac{|\vm|^2}{\vr} + 
\frac{1}{\ep^2} \vr e(\vr, S),
\]
whereas the relative energy can be written as \emph{Bregmann distance}
\begin{align}
E_\ep &\left( \vr, S, \vm \Big| \tvr , \tS, \tvm \right) = E_\ep(\vr, S, \vm) - \left< \partial_{\vr, S, \vm} E_\ep(\tvr, \tS, \tvm) ; (\vr - \tvr, S - \tS, \vm - \tvm) \right> - E_\ep(\tvr, \tS, \tvm)
\nonumber
\end{align}
associated to the convex functional $(\vr, S, \vm) \mapsto E_\ep (\vr, S, \vm)$, see \cite[Chapter~3, Section~3.1]{FeiNovOpen}.
Finally, as stated in \cite{ChauFei}, any weak solution in the sense of Definition \ref{DL1} satisfies the relative \textcolor{black}{energy} inequality 
 \begin{align}
	&\left[ \intO{ E_\ep \left(\vr, \vt, \vu \Big| \tvr, \tvt, \tvu \right) } \right]_{t = 0}^{t = \tau} \br 
	&+ \int_0^\tau \intO{ \frac{\tvt}{\vt} \left( \mathbb{S} (\vt, \Ds \vu) : \Ds \vu + \frac{1}{\ep^2} \frac{\kappa(\vt) |\Grad \vt|^2 }{\vt} \right) } \dt \br 
	&\leq - \frac{1}{\ep^2} \int_0^\tau \intO{ \left( \vr (s - s(\tvr, \tvt)) \partial_t \tvt + \vr (s - s(\tvr, \tvt)) \vu \cdot \Grad \tvt -
		\left( \frac{\kappa (\vt) \Grad \vt}{\vt} \right) \cdot \Grad \tvt \right) } \dt \br 
	&- \int_0^\tau \intO{ \Big[ \vr (\vu - \tvu) \otimes (\vu - \tvu) + \frac{1}{\ep^2} p(\vr, \vt) \mathbb{I} - \mathbb{S}(\vt, \Ds \vu) \Big] : \Ds \tvu } \dt \br 
	&+ \int_0^\tau \intO{ \vr \left[ \frac{1}{\ep} \Grad G  - \partial_t \tvu - (\tvu \cdot \Grad) \tvu  \right] \cdot (\vu - \tvu) } \dt \br 
	&+ \frac{1}{\ep^2} \int_0^\tau \intO{ \left[ \left( 1 - \frac{\vr}{\tvr} \right) \partial_t p(\tvr, \tvt) - \frac{\vr}{\tvr} \vu \cdot \Grad p(\tvr, \tvt) \right] } \dt
	\label{L4}
\end{align}
for a.a.~$\tau > 0$ and any trio of continuously differentiable functions $(\tvr, \tvt, \tvu)$ satisfying
\begin{equation} \label{L5}
	\tvr > 0,\ \tvt > 0,\ \tvt|_{\partial \Omega} = \vtB, \ \tvu|_{\partial \Omega} = 0.
\end{equation}
Here, we have denoted $\Ds \vu=\frac12(\Grad \vu +\Grad^t \vu)$.

\subsection{Constitutive relations}
\label{EOS}

The existence theory developed in \cite{ChauFei}, \cite{FeiNovOpen} is conditioned by certain restrictions imposed on the 
constitutive relations (state equations) similar to those introduced in the monograph
\cite[Chapters~1,2]{FeNo6A}. Specifically, the \emph{equation of state} takes the form
\begin{equation} \label{PES}
p(\vr, \vt) = p_{\rm m} (\vr, \vt) + p_{\rm rad}(\vt), 
\end{equation}
where $p_{\rm m}$ is the pressure of a general \emph{monoatomic} gas, 
\begin{equation} \label{con1}
	p_{\rm m} (\vr, \vt) = \frac{2}{3} \vr e_{\rm m}(\vr, \vt),
\end{equation}
enhanced by the radiation pressure 
\[
p_{\rm rad}(\vt) = \frac{a}{3} \vt^4,\ a > 0.
\]
Accordingly, the internal energy reads 
\[
e(\vr, \vt) = e_{\rm m}(\vr, \vt) + e_{\rm rad}(\vr, \vt),\ e_{\rm rad}(\vr, \vt) = \frac{a}{\vr} \vt^4.
\]
Moreover, using several physical principles it was shown in \cite[Chapter~1]{FeNo6A}:

\begin{itemize}
	
	\item \emph{Gibbs' relation} together with \eqref{con1} yield 
	\[
	p_{\rm m} (\vr, \vt) = \vt^{\frac{5}{2}} P \left( \frac{\vr}{\vt^{\frac{3}{2}}  } \right)
	\]
	for a certain $P \in C^1[0,\infty)$.
	Consequently, 
	\begin{equation} \label{w9}
		p(\vr, \vt) = \vt^{\frac{5}{2}} P \left( \frac{\vr}{\vt^{\frac{3}{2}}  } \right) + \frac{a}{3} \vt^4,\ 
		e(\vr, \vt) = \frac{3}{2} \frac{\vt^{\frac{5}{2}} }{\vr} P \left( \frac{\vr}{\vt^{\frac{3}{2}}  } \right) + \frac{a}{\vr} \vt^4, \ a > 0.
	\end{equation}
	
	\item \emph{Hypothesis of thermodynamic stability} \eqref{HTS} expressed in terms of 
	$P$ gives rise to
	\begin{equation} \label{w10}
		P(0) = 0,\ P'(Z) > 0 \ \mbox{for}\ Z \geq 0,\ 0 < \frac{ \frac{5}{3} P(Z) - P'(Z) Z }{Z} \leq C \ \mbox{for}\ Z > 0.
	\end{equation} 	
	In particular, the function $Z \mapsto P(Z)/ Z^{\frac{5}{3}}$ is decreasing, and we suppose 
	\begin{equation} \label{w11}
		\lim_{Z \to \infty} \frac{ P(Z) }{Z^{\frac{5}{3}}} = p_\infty > 0.
	\end{equation}
	
	\item 
	Accordingly, the associated entropy
	takes the form 
	\begin{equation} \label{w12}
		s(\vr, \vt) = s_{\rm m}(\vr, \vt) + s_{\rm rad}(\vr, \vt),\ s_{\rm m} (\vr, \vt) = \mathcal{S} \left( \frac{\vr}{\vt^{\frac{3}{2}} } \right),\ s_{\rm rad}(\vr, \vt) = \frac{4a}{3} \frac{\vt^3}{\vr}, 
	\end{equation}
	where 
	\begin{equation} \label{w13}
		\mathcal{S}'(Z) = -\frac{3}{2} \frac{ \frac{5}{3} P(Z) - P'(Z) Z }{Z^2} < 0.
	\end{equation}
\item	In addition, we impose the \emph{Third law of thermodynamics}, cf. Belgiorno \cite{BEL1}, \cite{BEL2}, requiring the entropy to vanish 
	when the absolute temperature approaches zero, 
	\begin{equation} \label{w14}
		\lim_{Z \to \infty} \mathcal{S}(Z) = 0.
	\end{equation}
	
\end{itemize}

Finally, we
suppose the transport coefficients are continuously differentiable functions satisfying
\begin{align} 
	0 < \underline{\mu}(1 + \vt) &\leq \mu(\vt),\ |\mu'(\vt)| \leq \Ov{\mu}, \br 
	0 \leq \eta (\vt) &\leq \Ov{\eta}(1 + \vt), \br
	0 < \underline{\kappa} (1 + \vt^\beta) &\leq \kappa (\vt) \leq \Ov{\kappa}(1 + \vt^\beta), 
	\ \mbox{where}\ \beta > 6. \label{w16}
\end{align}

As a consequence of the above hypotheses, we get the following estimates:
\begin{align} 
	\vr^{\frac{5}{3}} + \vt^4 \aleq \vr e(\vr, \vt) &\aleq 	1+ \vr^{\frac{5}{3}} + \vt^4, \label{L5b} \\
	s_{\rm m}(\vr, \vt) &\aleq \left( 1 + |\log(\vr)| + [\log(\vt)]^+ \right), \label{L5a}
\end{align} 
see \cite[Chapter~3, Section~3.2]{FeNo6A}. Here and hereafter, the symbol $a \aleq b$ means there is a positive constant $C> 0$ such that $a \le C b$.

We report the existence result proved in \cite[Theorem 4.2]{ChauFei}.

\begin{Proposition}[{\bf Primitive system, global existence}] \label{wP1}
	
	Let the thermodynamic functions $p$, $e$, $s$ and the transport coefficients $\mu$, $\eta$, $\kappa$ satisfy the hypotheses \eqref{PES}--\eqref{w16}.
	Let $\vt_B \in C^2(\partial \Omega)$. Suppose the initial data $\vr_0$, $\vt_0$, $\vu_0$ satisfy 
	\[
	\vr_0 > 0,\ \vt_0 > 0,\ \intO{ E_\ep(\vr_0, \vt_0, \vu_0 ) } < \infty.
	\]
	
	Then for any $T > 0$, the Navier--Stokes--Fourier system \eqref{i1}--\eqref{i7} admits a weak solution 
	$(\vr, \vt, \vu)$ in $(0,T) \times \Omega$ in the sense specified in Definition \ref{DL1}.

	\end{Proposition}
	
	\begin{Remark} \label{R1a}
	
	{\tc
		The constitutive restrictions imposed through \eqref{PES}--\eqref{w16}  are necessary for the global in time existence result stated in Proposition \ref{wP1}. They are of technical nature, and, in our opinion, 
	without significant impact of the asymptotic limit stated below. In particular, the relations \eqref{PES}, \eqref{con1} can be replaced by more general 
	ones specified in the monograph \cite[Chapter~1, Sections~1.4.2, 1.4.3]{FeNo6A}. The radiation pressure component is essential only for the existence theory eliminating hypothetical temperature 
	oscillations in the vacuum zones. The fact that its presence has no influence on the form of the asymptotic problem can be demonstrated by considering
	an extra scaling $a = a(\ep) \to 0$ as $\ep \to 0$. 
	}

	\end{Remark}

\section{Strong solutions to the target system} 
\label{OB}

Our analysis requires the existence of \emph{regular} solutions to the Oberbeck--Boussinesq system \eqref{ObBs1}--\eqref{i14B}. The relevant result was proved in \cite[Theorem 2.3, Theorem 3.1]{AbbFei22}. 

	\begin{Proposition}[\bf Strong solutions to target system] \label{OBP1}
	
	Suppose that 
	\begin{equation} \label{ss1}
		G \in W^{1, \infty}({\Omega}),\ \Theta_B \in C^2(\Ov{\Omega}),	
	\end{equation}
	and  
	\begin{align} 
		&\Theta_0 \in W^{2,p}(\Omega),\ \vU_0 \in W^{2,p}(\Omega; \mathbb{R}^d),\ \Div \vU_0 = 0, 
		\ \mbox{for any}\ 1 \leq p < \infty, \br
		&\mbox{together with the compatibility conditions} \br
		&\vU_0 = 0,\ \Theta_0 + \frac{\lambda}{1 - \lambda} \avintO{ \Theta_0 } = 
		\Theta_B \ \mbox{on}\ \partial \Omega.
		\label{M3}
	\end{align}
	
	Then there exists $T_{\rm max} > 0$, $T_{\rm max} = \infty$ if $d = 2$, such that
 the OB system \eqref{ObBs1}--\eqref{i14B} with  the initial data
	\[
	\vU(0, \cdot) = \vU_0,\ \Theta(0, \cdot) = \Theta_0,
	\]
	admits a strong solution
	  $\vU$, $\Theta$ in the regularity class
	\begin{align} 
		\vU \in L^p(0,T; W^{2,p}(\Omega; \mathbb{R}^d)),\ \partial_t \vU \in L^p(0,T; L^{p}(\Omega; \mathbb{R}^d)),\ \Pi \in L^p(0,T; W^{1,p}(\Omega)), \br 
		\Theta \in L^p(0,T; W^{2,p}(\Omega)),\ 
		\partial_t \Theta \in L^p(0,T; L^{p}(\Omega; \mathbb{R}^d))
		\label{M4}	
	\end{align}	
	for any $1 \leq p < \infty$ and any $0 < T < T_{\rm max}$.	
	
\end{Proposition}

\begin{Remark} \label{regul}
	
Strictly speaking, the existence result in \cite{AbbFei22} requires, in addition to \eqref{ss1},
\begin{equation} \label{poten}
\Del G = 0.
\end{equation}
On the one hand, in applications, $G$ represents the gravitational potential therefore \eqref{poten} is automatically satisfied. On the other hand, the proof presented in \cite{AbbFei22} can be easily modified to accommodate the general case.
	
	\end{Remark}

Given the parabolic character of the OB system, the strong solutions are in fact classical if higher regularity of the data is required, cf. \cite[Theorem 4.1]{AbbFei22}. 

\section{Asymptotic limit, main result}
\label{M}

We are ready to state our main result concerning the singular limit $\ep \to 0$ in the primitive NSF system.

\begin{mdframed}[style=MyFrame]
	
	\begin{Theorem}[{\bf Singular limit}] \label{MT1}
		
		Let the constitutive hypotheses \eqref{PES}--\eqref{L5a} be satisfied. Let the data belong to the regularity class 
		\[
		G  \in W^{1,\infty}(\Omega), \ \Theta_B \in C^2(\Ov{\Omega}).
		\]
		Let $(\vre, \vte, \vue)_{\ep > 0}$ be a family of weak solutions to the NSF system \eqref{i1}--\eqref{i7}, with the boundary data 
		\begin{equation} \label{LV6a}
		\vue|_{\partial \Omega} = 0,\ 
		\vte|_{\partial \Omega} = \Ov{\vt} + \ep \Theta_B,\ \Ov{\vt} > 0 \ \mbox{constant,}
		\end{equation}
		and the initial data
		\begin{equation} \label{illP}
		\vr_{\ep}(0, \cdot) = \Ov{\vr} + \ep \vr_{0, \ep},\ 
		 \Ov{\vr} > 0 \ \mbox{constant},\ \intO{ \vr_{0, \ep} } = 0,\ \vte(0, \cdot) = \Ov{\vt} + \ep \vt_{0, \ep},\ \vue(0, \cdot) = \vu_{0, \ep},
		\end{equation}
		where, in addition, 
		\begin{align}
		\| \vr_{0, \ep} \|_{L^\infty(\Omega) } &\aleq 1,\ \vr_{0, \ep} \to r_0 \ \mbox{in}\ L^1(\Omega), \br
		\| \vt_{0, \ep} \|_{L^\infty(\Omega) } &\aleq 1,\ \vt_{0, \ep} \to \mathfrak{T}_0 \ \mbox{in}\ L^1(\Omega), \br 
		\| \vu_{0, \ep} \|_{L^\infty(\Omega; \mathbb{R}^d) } &\aleq 1,\ \vu_{0, \ep} \to \vU_0 \ \mbox{in}\ L^1(\Omega; \mathbb{R}^d), 
		 \label{M5}
		\end{align}
	and
	\begin{align} 
		\mathfrak{T}_0 &\in W^{2,p}(\Omega),\ \vU_0 \in W^{2,p}(\Omega; \mathbb{R}^d), 
		\ \mbox{for any}\ 1 \leq p < \infty, \ \Div \vU_0 = 0,\br
		\vU_0 &= 0,\ \mathfrak{T}_0 = 
		\Theta_B \ \mbox{on}\ \partial \Omega,
		\label{M6}
	\end{align}	
\begin{equation} \label{M6a}
\frac{\partial p(\Ov{\vr}, \Ov{\vt} ) }{\partial \vr} \Grad r_0 + 
\frac{\partial p(\Ov{\vr}, \Ov{\vt} ) }{\partial \vt} \Grad \mathfrak{T}_0	 = \Ov{\vr} \Grad G.
\end{equation}

Then 
\begin{align} 
	\frac{\vre - \Ov{\vr}}{\ep} &\to r \ \mbox{in}\ L^\infty(0,T; L^1(\Omega)), \br
		\frac{\vte - \Ov{\vt}}{\ep} &\to \mathfrak{T} \ \mbox{in}\ L^\infty(0,T; L^1(\Omega)), \br
		\sqrt{\vre} \vue &\to \sqrt{\Ov{\vr}} \vc{U} \ \mbox{in}\  L^\infty(0,T; L^2(\Omega; \mathbb{R}^d)),
		\label{M7}
	\end{align}
as $\ep \to 0$ for any $T < T_{\rm max}$, where
\[
r,\ \Theta = \mathfrak{T} - \lambda(\Ov{\vr}, \Ov{\vt}) \avintO{ \mathfrak{T} },\ \vU	
\]
is the strong solution of the OB system \eqref{ObBs1}--\eqref{i14B} with the initial data 
$\vU(0, \cdot) = \vU_0$, $\Theta_0 = \mathfrak{T}_0 - \lambda \avintO{ \mathfrak{T}_0 }$ defined on $[0, T_{\rm max})$.
		\end{Theorem}

	\end{mdframed}

Note that hypotheses {\cred \eqref{M5}--\eqref{M6a}} correspond to \emph{well prepared} data.
The rest of the paper is devoted to the proof of Theorem \ref{MT1}.

\section{Basic energy estimates}
\label{B}

In order to perform the limit claimed in Theorem \ref{MT1}, we need bounds on the sequence 
$(\vre, \vte, \vue)_{\ep > 0}$ independent of $\ep$. We start by introducing the notation borrowed from  
\cite{FeNo6A} distinguishing the ``essential'' and ``residual'' range of the 
thermostatic variables $(\vr, \vt)$. Specifically, given a compact set 
\begin{equation*}
K \subset \left\{ (\vr, \vt) \in \R^2 \ \Big| \ \vr > 0, \vt > 0 \right\}
\end{equation*}
and $\ep > 0$, we introduce 
\begin{equation*}
g_{\rm ess} = g \mathds{1}_{(\vre, \vte) \in K},\ 
g_{\rm res} = g - g_{\rm ess} = g \mathds{1}_{(\vre, \vte) \in \R^2 \setminus K}
\end{equation*}
for any measurable $g = g(t,x)$. Note carefully that this decomposition depends on $\ep$. As a matter of fact, the characteristic function $\mathds{1}_{(\vre, \vte) \in K}$ can be replaced by its smooth regularization by a suitable convolution kernel. 

In the subsequent analysis, we consider
\[
K = \Ov{\mathcal{U}(\Ov{\vr}, \Ov{\vt})} \subset (0, \infty)^2,\ {\cred \mathcal{U} (\Ov{\vr}, \Ov{\vt} )} \ \mbox{- an open neighborhood of}\ (\Ov{\vr}, \Ov{\vt}).
\]
As shown in \cite[Chapter~5, Lemma~5.1]{FeNo6A}, the relative energy enjoys the following coercivity properties:
\[
E_{\ep} \left( \vr, \vt, \vu \Big| \tvr, \tvt, \tvu \right) \geq 
C \left( \frac{ |\vr - \tvr|^2 }{\ep^2} + \frac{ |\vt - \tvt|^2 }{\ep^2} + |\vu - \tvu |^2 \right)
\]
if $(\vr, \vt) \in K = \Ov{\mathcal{U}(\Ov{\vr}, \Ov{\vt})}$, $(\tvr, \tvt) \in \mathcal{U}(\Ov{\vr}, \Ov{\vt})$,
\[
E_{\ep} \left( \vr, \vt, \vu \Big| \tvr, \tvt, \tvu \right) \geq 
C \left( \frac{1}{\ep^2} + \frac{1}{\ep^2} \vr e(\vr, \vt) + \frac{1}{\ep^2} \vr |s(\vr, \vt)| + \vr |\vu|^2 \right)
\]
whenever $(\vr, \vt) \in R^2 \setminus \Ov{\mathcal{U}(\Ov{\vr}, \Ov{\vt})}$, $(\tvr, \tvt) \in \mathcal{U}(\Ov{\vr}, \Ov{\vt})$. The constant $C$ depends on $K$ and the distance 
\[
	\sup_{t,x} {\rm dist} \left[ (\tvr (t,x), \tvt (t,x) ) ; \partial K \right]. 	\]
	
In other words,
\begin{align} 
E_{\ep} \left( \vre, \vte, \vue \Big| \tvr, \tvt, \tvu \right)_{\rm ess} &\geq 
C \left( \frac{ |\vre - \tvr|^2 }{\ep^2} + \frac{ |\vte - \tvt|^2 }{\ep^2} + |\vue - \tvu |^2 \right)_{\rm ess}, 
	\label{BB1} \\ 
E_{\ep} \left( \vre, \vte, \vue \Big| \tvr, \tvt, \tvu \right)_{\rm res} &\geq 
C  \left( \frac{1}{\ep^2} + \frac{1}{\ep^2} \vre e(\vre, \vte) + \frac{1}{\ep^2} \vre |s(\vre, \vte)| + \vre |\vue|^2 \right)_{\rm res}. 
\label{BB2}	
	\end{align}

\subsection{Energy estimates}

In agreement with hypothesis \eqref{illP}, we have 
\begin{equation} \label{L1}
	\intO{ E_\ep \left( \vre(0, \cdot), \vte(0, \cdot), \vue(0, \cdot) \Big| \Ov{\vr}, \Ov{\vt} + \ep \Theta_B , 0 \right) } \aleq 1 \ \mbox{independently of}\ \ep \to 0,
		\end{equation}
where $\Theta_B = \Theta_B(x)$ is a suitable extension of the temperature boundary condition inside $\Omega$. 
Plugging this ansatz in the relative energy inequality \eqref{L4} we obtain
 \begin{align}
	&\left[ \intO{ E_\ep \left(\vre, \vte, \vue \Big| \Ov{\vr}, \Ov{\vt} + \ep \Theta_B, 0 \right) } \right]_{t = 0}^{t = \tau} \br 
	&+ \int_0^\tau \intO{ \frac{\Ov{\vt} + \ep \Theta_B}{\vte} \left( \mathbb{S} (\vte, \Ds \vue) : \Ds \vue + \frac{1}{\ep^2} \frac{\kappa (\vte) |\Grad \vte |^2 }{\vte} \right) } \dt \br 
	&\leq - \frac{1}{\ep} \int_0^\tau \intO{ \left(  \vre (s(\vre, \vte) - s(\Ov{\vr}, \Ov{\vt} + \ep{\Theta_B})) \vue \cdot \Grad \Theta_B - \frac{\kappa(\vte) \Grad \vte }{\vte} \cdot \Grad \Theta_B
 \right) } \dt \br 
	&\quad + \int_0^\tau \intO{ \vre \frac{1}{\ep} \Grad G  \cdot \vue  } \dt \br 
	&\quad {\cred - \frac{1}{\ep} }\int_0^\tau \intO{ \left( \frac{ \partial p(\Ov{\vr}, \Ov{\vt} + \ep \Theta_B )}{\partial \vt} - \frac{ \partial p(\Ov{\vr}, \Ov{\vt}  )}{\partial \vt} \right) \frac{\vre}{\Ov{\vr}} \vue \cdot \Grad \Theta_B  } \dt \br 
	&\quad - \frac{1}{\ep} \int_0^\tau \intO{  \frac{ \partial p(\Ov{\vr}, \Ov{\vt}  )}{\partial \vt} \frac{\vre}{\Ov{\vr}} \vue \cdot \Grad \Theta_B  } \dt.
	\label{L2}
\end{align}

Our goal is to control the integrals on the right--hand side to apply Gr\" onwall's argument. To this end, we fix the compact set $K$ determining the essential and residual component to contain the point $(\Ov{\vr}, \Ov{\vt})$ in its interior.
In particular, the same is true for the range of the function $(\Ov{\vr}, \Ov{\vt} + \ep \Theta_B)$ as soon as $\ep > 0$ is small enough. Accordingly, we will systematically use the coercivity of the relative energy 
$E_\ep$ stated in \eqref{BB1}, \eqref{BB2} in the estimates below. In particular, we have the 
	estimate
\begin{align}
\left| \left[ B(\vre, \vte) - B(\Ov{\vr}, \Ov{\vt} + \ep \Theta_B) \right]_{\rm ess} \right| 
&\aleq \left[ E_\ep \left(\vre, \vte, \vue \Big| \Ov{\vr}, \Ov{\vt} + \ep \Theta_B, 0 \right) \right]_{\rm ess} 
\br &\leq E_\ep \left(\vre, \vte, \vue \Big| \Ov{\vr}, \Ov{\vt} + \ep \Theta_B, 0 \right)
\nonumber
\end{align}
for any $B = B(\vr, \vt)$ locally Lipschitz in $(0, \infty)^2$.

\subsubsection{Estimates}

{\bf Step 1:}  First, 
\begin{align}
\frac{1}{\ep} &\intO{ \left|  \vre (s(\vre, \vte) - s(\Ov{\vr}, \Ov{\vt} + \ep{\Theta_B})) \vue \cdot \Grad \Theta_B \right| }	\br 
&\aleq \frac{1}{\ep} \intO{ \left| \left[ \vre (s(\vre, \vte) - s(\Ov{\vr}, \Ov{\vt} + \ep{\Theta_B})) \vue \right]_{\rm ess}  \right| } \br &\quad + 
\frac{1}{\ep} \intO{ \left| \left[ \vre (s(\vre, \vte) - s(\Ov{\vr}, \Ov{\vt} + \ep{\Theta_B})) \vue  \right]_{\rm res}  \right| },
\nonumber
	\end{align}
where 
\begin{align}
\frac{1}{\ep} &\intO{ \left| \left[ \vre (s(\vre, \vte) - s(\Ov{\vr}, \Ov{\vt} + \ep{\Theta_B})) \vue \right]_{\rm ess}  \right| } \aleq	
\frac{1}{\ep^2} \intO{ \left| \left[ (s(\vre, \vte) - s(\Ov{\vr}, \Ov{\vt} + \ep{\Theta_B})) \right]_{\rm ess} \right|^2  } \br 
&+ \intO{ \vre |\vue|^2 } \aleq \intO{ E_\ep \left( \vre, \vte, \vue \Big| \Ov{\vr}, \Ov{\vt} + \ep \Theta_B, 0 \right) },
	\label{L3}
	\end{align}
and 
\begin{align} 
\frac{1}{\ep} &\intO{ \left| \left[ \vre (s(\vre, \vte) - s(\Ov{\vr}, \Ov{\vt} + \ep{\Theta_B})) \vue \right]_{\rm res}  \right| }\br &\aleq 
\frac{1}{\ep} \intO{ \left[ \vre  |\vue| \right]_{\rm res} } + 
\frac{1}{\ep} \intO{ \left[ \vre s_{\rm m}(\vre, \vte) |\vue| \right]_{\rm res} } + 
\frac{1}{\ep} \intO{ \left[  \vte^3 |\vue| \right]_{\rm res} }. 
\nonumber	
	\end{align}
Furthermore, 
\begin{equation} \label{L4a}
	\frac{1}{\ep} \intO{ \left[ \vre  |\vue| \right]_{\rm res} } \aleq \frac{1}{\ep^2} \intO{ [\vr_\ep]_{\rm res} } + \intO{\vre |\vue|^2} \aleq
 \intO{ E_\ep \left( \vre, \vte, \vue \Big| \Ov{\vr}, \Ov{\vt} + \ep \Theta_B, 0 \right) }
\end{equation}
In view of the bounds \eqref{L5b}, \eqref{L5a}, 
\begin{align} 
\frac{1}{\ep} \intO{ \left[ \vre s_{\rm m}(\vre, \vte) |\vue| \right]_{\rm res} } &\aleq \frac{1}{\ep^2} \intO{ \left[ \vre s^2_{\rm m}(\vre, \vte) \right]_{\rm res} } + \intO{\vre |\vue|^2}  \br
&\aleq
\intO{ E_\ep \left( \vre, \vte, \vue \Big| \Ov{\vr}, \Ov{\vt} + \ep \Theta_B, 0 \right) }.
\label{b1}	
	\end{align}

Finally, 
\begin{align} 
\frac{1}{\ep} \intO{ \left[  \vte^3 |\vue| \right]_{\rm res} } \aleq \delta \| \vue \|_{W^{1,2}(\Omega; \mathbb{R}^d)}^2 + \frac{  C(\delta) }{\ep^2} \intO{ [\vte^6]_{\rm res} } 
\nonumber	
	\end{align}
for any $\delta > 0$. Thus if $\delta > 0$ is chosen small enough, the first integral is controlled by the viscosity dissipation on the left--hand side of \eqref{L2}. Next, in accordance with hypothesis \eqref{w16},
\begin{equation} \label{b2}
\intO{ \frac{\kappa (\vte) |\Grad \vte|^2 }{\vte^2} } \ageq 
\intO{ |\Grad \log(\vte) |^2 + |\Grad \vte^{\frac{\beta}{2}} |^2 },\ \beta > 6.
\end{equation}
Consequently, as the measure of the residual set is controlled by the relative energy (cf. \eqref{BB2}), we get 
\begin{align} 
\frac{ 1 }{\ep^2} &\intO{ [\vte^6]_{\rm res} } \aleq \frac{1}{\ep^2} \intO{ [\vte^3 - \Ov{\vt}^3 ]_{\rm res}^2 } + \frac{1}{\ep^2} \intO{ [\Ov{\vt}^3]_{\rm res}^2 } \br
&\aleq
\frac{1}{\ep^2}  \intO{ |\Grad \vte^{3}  |^2 } + \frac{1}{\ep^2} \int_{\Omega} [\vte^3 - \Ov{\vt}^3 ]_{\rm ess}^2 + \frac{1}{\ep^2} \intO{ [\Ov{\vt}^3]_{\rm res}^2 }\br & \aleq 
\frac{\delta}{\ep^2} \intO{ \frac{\kappa (\vte) |\Grad \vte|^2 }{\vte^2} } + C(\delta)\intO{ E_\ep \left( \vre, \vte, \vue \Big| \Ov{\vr}, \Ov{\vt} + \ep \Theta_B, 0 \right) }  
\label{b4} 
\end{align}
for any $\delta > 0$.

\bigskip

{\bf Step 2:} In accordance with hypothesis \eqref{w16}, 
\[
\frac{1}{\ep} \left| \intO{ \frac{\kappa(\vte) }{\vte} \Grad \vte } \right| \aleq \frac{1}{\ep} \intO{ |\Grad (\log(\vte))| + \vte^{\beta - 1} |\Grad \vte| },
\]
where
\begin{align}
\frac{1}{\ep} \intO{ \left|  \Grad (\log(\vte)) \right| } \aleq \frac{\delta}{\ep^2} \intO{ \left|  \Grad (\log(\vte)) \right|^2 } + C(\delta)
\label{b5}
\end{align}
for any $\delta> 0$; hence the integral is controlled by dissipation. 

Next, 
\begin{align}
\frac{1}{\ep} \intO{ \vte^{\beta - 1}{\Grad \vte} } = 
\frac{1}{\ep} \intO{ \vt^{\frac{\beta}{2}} \Grad \vte^{\frac{\beta}{2}} } \leq \frac{\delta}{\ep^2} \intO{ |\Grad \vte^{\frac{\beta}{2}} |^2 } + C(\delta) \intO{ |\vte^{\frac{\beta}{2} }|^2 }, 	
\label{b6}
	\end{align}
where the first term is controlled by dissipation and the second one by Poincar\' e's inequality 
\begin{align}
\intO{ |\vte^{\frac{\beta}{2} } |^2 } \aleq \intO{ |\Grad \vte^{\frac{\beta}{2}} |^2 } + \int_{\partial \Omega} ( \Ov{\vt} + \ep \Theta_B )^\beta \ \D \sigma_x.	
	\label{b7}
	\end{align}

\bigskip

{\bf Step 3:} We have 
\begin{align}
\frac{1}{\ep} \int_0^\tau \intO{ \vre \Grad G \cdot \vue } \dt &= - \frac{1}{\ep} \int_0^\tau \intO{ G \Div (\vre \vue) } \dt = 
{\cred \int_0^\tau \partial_t  \intO{ \frac{1}{\ep} (\vre - \Ov{\vr}) G } \dt }\br  = & \left[ \intO{ \frac{\vre - \Ov{\vr}}{\ep}  G } \right]_{t = 0}^{t=\tau}.
\nonumber
\end{align}
Seeing that 
\begin{align}
E_\ep \left( \vre, \vte, \vue \Big| \Ov{\vr}, \Ov{\vt} + \ep \Theta_B, 0 \right) - c_1 &\aleq 
E_\ep \left( \vre, \vte, \vue \Big| \Ov{\vr}, \Ov{\vt} + \ep \Theta_B, 0 \right)
- \intO{ \frac{\vre - \Ov{\vr}}{\ep}  G } \br &\aleq E_\ep \left( \vre, \vte, \vue \Big| \Ov{\vr}, \Ov{\vt} + \ep \Theta_B, 0 \right) + 1
\end{align}
we can add this term to the relative energy on the left--hand side of \eqref{L4}.
\bigskip 

\noindent
{\bf Step 4:}
\begin{align}
\frac{1}{\ep} &\left| \intO{ \left( \frac{ \partial p(\Ov{\vr}, \Ov{\vt} + \ep \Theta_B )}{\partial \vt} - \frac{ \partial p(\Ov{\vr}, \Ov{\vt}  )}{\partial \vt} \right) \frac{\vre}{\Ov{\vr}} \vue \cdot \Grad \Theta_B  } \right| \br 
&\aleq \intO{ \vre |\vue| } \aleq \intO{ \vre } + \intO{ \vre |\vue|^2 } \aleq \intO{ E_\ep \left( \vre, \vte, \vue \Big| \Ov{\vr}, \Ov{\vt} + \ep \Theta_B, 0 \right) } + 1.
\label{b8}
\end{align}

\bigskip 

\noindent
{\bf Step 5:} The last integral on the right--hand side of \eqref{L2} can be handled exactly as in Step 3.

\subsubsection{Conclusion, uniform bounds}
\label{cbe}

In view of the estimates obtained in the previous section, we may apply Gr\" onwall's lemma to the relative energy inequality \eqref{L2}. As the initial data  satisfy \eqref{L1}, we deduce
the following bounds independent of the scaling parameter $\ep \to 0$: 

\begin{align}
	{\rm ess} \sup_{t \in (0,T)} \intO{ E_\ep \left( \vre, \vte, \vue \Big| \Ov{\vr}, \Ov{\vt} + \ep \Theta_B, 0 \right) } &\aleq 1, \label{be1} \\
	\int_0^T \| \vue \|^2_{W^{1,2}_0 (\Omega; \mathbb{R}^d) } \dt &\aleq 1, \label{be2} \\
	\frac{1}{\ep^2} \int_0^T \left( \| \Grad \log (\vte) \|^2_{L^2(\Omega; \mathbb{R}^d)} + \| \Grad \vte^{\frac{\beta}{2}} \|^2_{L^2(\Omega; \mathbb{R}^d)} \right) &\aleq 1.
	\label{be3}
	\end{align}

Next, it follows from \eqref{be1} that the measure of the residual set shrinks to zero, specifically
\begin{equation} \label{be4}
	\frac{1}{\ep^2} {\rm ess} \sup_{t \in (0,T)} \intO{ [1]_{\rm res} } \aleq 1.
\end{equation}

In addition, we get from \eqref{be1}:
\begin{align}
	{\rm ess} \sup_{t \in (0,T)} \intO{ \vre |\vue|^2 } &\aleq 1, \br 
	{\rm ess} \sup_{t \in (0,T)} \left\| \left[ \frac{\vre - \Ov{\vr}}{\ep} \right]_{\rm ess} \right\|_{L^2(\Omega)} &\aleq 1, \br
	{\rm ess} \sup_{t \in (0,T)} \left\| \left[ \frac{\vte - \Ov{\vt}}{\ep} \right]_{\rm ess} \right\|_{L^2(\Omega)} &\aleq 1, \br
	\frac{1}{\ep^2} {\rm ess} \sup_{t \in (0,T)} \| [\vre]_{\rm res} \|^{\frac{5}{3}}_{L^{\frac{5}{3}}(\Omega)} + \frac{1}{\ep^2}  {\rm ess} \sup_{t \in (0,T)} \| [\vte]_{\rm res} \|^{4}_{L^{4}(\Omega)}	 &\aleq 1.
\label{be5}	
	\end{align}

Combining \eqref{be3}, \eqref{be4}, and \eqref{be5}, we conclude 
\begin{equation} \label{be6}
	\int_0^T \left\| \frac{\log(\vte) - \log(\Ov{\vt})}{\ep} \right\|^2_{W^{1,2}(\Omega)} \dt + \int_0^T \left\| \frac{ \vte - \Ov{\vt} }{\ep} \right\|^2_{W^{1,2}(\Omega)} \dt \aleq 1.
	\end{equation}

Finally, we claim the bound on the entropy flux
\begin{equation} \label{be7} 
\int_0^T \left\| \left[ \frac{\kappa (\vte) }{\vte} \right]_{\rm res} \frac{\Grad \vte }{\ep} \right\|^q_{L^q(\Omega; \mathbb{R}^d)} \dt \aleq 1 \ \mbox{for some}\ q > 1.
\end{equation} 
Indeed we have 
\[
\tc
\left|  \left[ \frac{\kappa (\vte) }{\vte} \right]_{\rm res} \frac{\Grad \vte }{\ep} \right| \aleq  \frac{1}{\ep}
\left| \Grad \log (\vte) \right| + \frac{1}{\ep} \left| \left[ \vte^{\frac{\beta}{2}} \Grad \vte^{\frac{\beta}{2}} \right]_{\rm res} \right|, 
\]	
where the former term on the right--hand side is controlled via \eqref{be6}. As	for the latter, we deduce from 
\eqref{be3} that  
\[
\left\| \frac{1}{\ep} \Grad \vte^{\frac{\beta}{2}} \right\|_{L^2((0,T) \times \Omega; \mathbb{R}^d)} \aleq 1;
\]
hence it is enough to check  
\begin{equation} \label{be8}
\tc
\left\| \left[ \vte^{\frac{\beta}{2}} \right]_{\rm res}\right\|_{L^r ((0,T) \times \Omega)} \aleq 1 \ \mbox{for some}\ r > 2.
\end{equation}
To see \eqref{be8}, first observe that 
\[
\tc
{\rm ess} \sup_{t \in (0,T)} \| [\vte]_{\rm res} \|_{L^4(\Omega)} \aleq 1, 
\]
and,  in view of \eqref{be3} and Poincar\' e inequality, 
\[
\left\| 
\vte^{\frac{\beta}{2}} \right\|_{L^2(0,T; L^6(\Omega))} \aleq 1 \ \mbox{(for $d = 3$)}.
\]
Consequently, \eqref{be8} follows by interpolation.

\section{Convergence to the target system}

\label{WPL}

Our ultimate goal is to perform the limit $\ep \to 0$. We proceed in two steps. 

\subsection{Weak convergence} 

In view of the uniform bounds established in Section~\ref{cbe},  
we may infer 
\begin{align} 
	\vre &\to \Ov{\vr} \ \mbox{in}\ L^{\frac{5}{3}}(\Omega) \ \mbox{uniformly for}\ t \in (0,T),  \label{wc1} \\ 
	\vte &\to \Ov{\vt} \ \mbox{in}\ L^2(0,T; W^{1,2}(\Omega)), \label {wc2} \\
	\vue &\to \vu \ \mbox{weakly in}\ L^2(0,T; W^{1,2}_0(\Omega; \mathbb{R}^d)), \label{wc3}
\end{align}
where \eqref{wc3} may require extraction of a suitable subsequence. As we shall eventually see,  
the limit velocity $\vu = \vU$ is unique so that the convergence is, in fact, unconditional.
In addition, we may let $\ep \to 0$ in the weak formulation of the equation of continuity \eqref{Lw4} to deduce 
\begin{equation} \label{wc4}
	\Div \vu = 0. 
\end{equation}

Next, we use \eqref{be5}, \eqref{be6} to obtain ({\it a priori} for suitable subsequences), 
\begin{align}
	\frac{ \vre - \Ov{\vr} }{\ep} &= \left[ \frac{ \vre - \Ov{\vr} }{\ep} \right]_{\rm ess} + 
	\left[ \frac{ \vre - \Ov{\vr} }{\ep} \right]_{\rm res}, \br 
	\left[ \frac{ \vre - \Ov{\vr} }{\ep} \right]_{\rm ess} &\to \mathfrak{R} \ \mbox{weakly-(*) in}\ 
	L^\infty(0,T; L^2(\Omega)),\br
	\left[ \frac{ \vre - \Ov{\vr} }{\ep} \right]_{\rm res} &\to 0 \ \mbox{in}\ L^\infty(0,T; L^{\frac{5}{3}}(\Omega)),
	\label{wc5}	
\end{align}
\begin{equation}
	\frac{ \vte - \Ov{\vt} }{\ep} \to  \mathfrak{T} \ \mbox{weakly in}\  L^2(0,T; W^{1,2}(\Omega)) 
	\ \mbox{and weakly-(*) in}\ L^\infty(0,T; L^2(\Omega)).
	\label{wc6}	
\end{equation}
Moreover, in view of \eqref{LV6a}, 
\begin{equation} \label{tep}
	\mathfrak{T}|_{\partial \Omega} = \Theta_B.
\end{equation}

Finally, we perform the limit in the rescaled momentum equation \eqref{i2} to deduce 
\begin{equation} \label{wc7}
	\frac{\partial p(\Ov{\vr}, \Ov{\vt}) }{\partial \vr} \Grad \mathfrak{R} + 
	\frac{\partial p(\Ov{\vr}, \Ov{\vt}) }{\partial \vt} \Grad \mathfrak{T} = \Ov{\vr}  \Grad G
\end{equation}
in the sense of distributions. In particular, it follows from \eqref{wc7} that 
\begin{equation} \label{wc8}
	\mathfrak{R} \in L^2(0,T; W^{1,2}(\Omega)).
\end{equation}

\subsection{Strong convergence}

First, it is more convenient to rewrite the target OB system in terms of the variable 
\[
\mathcal{T}, \ \mbox{where} \ \mathcal{T} - \lambda (\Ov{\vr}, \Ov{\vt} ) \avintO{ \mathcal{T} } = \Theta. 
\]
Accordingly, we get 
\begin{align} 
	\Div \vU &= 0, \br	
	\Ov{\vr} \Big( \partial_t \vU + \vU \cdot \Grad \vU \Big) + \Grad \Pi &= \Div \mathbb{S}(\Ov{\vt}, \Grad \vU) +  r \Grad G, \br
	\Ov{\vr} c_p(\Ov{\vr}, \Ov{\vt} ) \left( \partial_t \MTC + \vU \cdot \Grad \MTC \right)	- 
	\Ov{\vr} \ \Ov{\vt} \alpha(\Ov{\vr}, \Ov{\vt} ) \vU \cdot \Grad G
	&= \kappa(\Ov{\vt}) \Del \MTC + \Ov{\vt} \alpha (\Ov{\vr}, \Ov{\vt}) \frac{\partial p (\Ov{\vr}, \Ov{\vt})}{\partial \vt}  
	\partial_t \avintO{ \MTC },	
	\label{OBs}
\end{align}
together with the Boussinesq relation
\begin{equation} \label{OB1}
	\frac{\partial p(\Ov{\vr}, \Ov{\vt} ) }{\partial \vr} \Grad r + 
	\frac{\partial p(\Ov{\vr}, \Ov{\vt} ) }{\partial \vt} \Grad \MTC = \Ov{\vr} \Grad G,\ 
	\intO{ r } = 0,
\end{equation} 
the boundary conditions 
\begin{equation} \label{bc}
	\vU|_{\partial \Omega} = 0,\ \MTC|_{\partial \Omega} = \Theta_B,
\end{equation}
and the initial conditions
\begin{equation} \label{OBini}
	\vU (0, \cdot) = \vU_0,\ \MTC(0, \cdot) = \mathfrak{T}_0.
	\end{equation}
In accordance with Proposition \ref{OBP1} and hypotheses \eqref{M6}, \eqref{M6a}, the problem 
\eqref{OBs}--\eqref{OBini} admits a unique regular solution on a time interval $[0, T_{\rm max})$, where 
$T_{\rm max} > 0$ and $T_{\rm max} = \infty$ if $d=2$.

\subsection{Relative energy}

To complete the proof of Theorem \ref{MT1}, we use the relative energy inequality \eqref{L4}, with the ansatz 
\[
E_\ep \left(\vre, \vte, \vue \ \Big| \Ov{\vr} + \ep r, \Ov{\vt} + \ep  \MTC, \vU \right). 
\]
In accordance with our choice of the initial data, 
\begin{equation} \label{idC}
\intO{ E_\ep \left(\vre, \vte, \vue \ \Big| \Ov{\vr} + \ep r, \Ov{\vt} + \ep \MTC, \vU \right) (0, \cdot) 
} \to 0 \ \mbox{as}\ \ep \to 0.
\end{equation}
{\cred Our goal is to show that $\displaystyle \intO{E_\ep \left(\vre, \vte, \vue \ \Big| \Ov{\vr} + \ep r, \Ov{\vt} + \ep \MTC, \vU \right)(\tau,\cdot)}\to 0$, which finally yields
\begin{align*}
\MTC=\lim_{\ep \to 0} \frac{\vte-\Ov{\vt}}{\ep}=\mathfrak{T},\ r=\lim_{\ep\to 0} \frac{\vre-\Ov{\vr}}{\ep}=\mathfrak{R},\ \lim_{\ep\to 0} \vue=\vU,
\end{align*}
and the trio $(\mathfrak{R},\mathfrak{T},\vU)=(r, \MTC, \vU)$ is the strong solution to the OB system \eqref{ObBs1}--\eqref{i14B}.
}

\noindent
{\bf Step 1:} Plugging our ansatz in the relative energy inequality \eqref{L4} and using 
$\Div \vU = 0$ we get  
 \begin{align}
	&\left[ \intO{ E_\ep \left(\vre, \vte, \vue \Big| \Ov{\vr} + \ep r, \Ov{\vt} + \ep \MTC, \vU \right) } \right]_{t = 0}^{t = \tau} \br 
	&+ \int_0^\tau \intO{ \frac{\Ov{\vt} + \ep \MTC}{\vte} \left( \mathbb{S} (\vte, \Ds \vue) : \Ds \vue + \frac{1}{\ep^2} \frac{\kappa (\vte) \Grad \vte \cdot \Grad \vte }{\vte} \right) } \dt \br 
	&\leq - \frac{1}{\ep} \int_0^\tau \intO{ \vre \Big[ (s(\vre, \vte) - s( \Ov{\vr} + \ep r,\Ov{\vt} + \ep \MTC ) \Big] \partial_t \MTC } \dt \br 
	&\quad - \frac{1}{\ep} \int_0^\tau \intO{ \vre \Big[ s(\vre, \vte) - s( \Ov{\vr} + \ep r, \Ov{\vt} + \ep \MTC) \Big] \vue \cdot \Grad \MTC  } \dt \br
	&\quad  + \frac{1}{\ep} \int_0^\tau \intO{ 
	 \frac{\kappa(\vte) }{\vte} \Grad \vte \cdot \Grad \MTC } \dt \br 
	&\quad - \int_0^\tau \intO{ \Big[ \vre (\vue - \vU) \otimes (\vue - \vU) - \mathbb{S}(\vte, \Ds \vue) \Big] : \Ds \vU } \dt \br 
	&\quad + \int_0^\tau \intO{ \vre \left[ \frac{1}{\ep} \Grad G  - \partial_t \vU - (\vU \cdot \Grad) \vU \right] \cdot (\vue - \vU) } \dt \br 
	&\quad + \frac{1}{\ep^2} \int_0^\tau \intO{ \left[ \left( 1 - \frac{\vre}{\Ov{\vr} + \ep r} \right) \partial_t p ( \Ov{\vr} + \ep r, \Ov{\vt} + \ep \MTC) - \frac{\vre}{\Ov{\vr} + \ep r} \vue \cdot \Grad p( \Ov{\vr} + \ep r, \Ov{\vt} + \ep \MTC) \right] } \dt.
	\label{L6}
\end{align}

\bigskip 
\noindent
{\bf Step 2:} As $r$, $\vU$ satisfy the momentum equation
\[
- \Ov{\vr} \left( \partial_t \vU + \vU \cdot \Grad \vU \right) = \Grad \Pi - \Div \mathbb{S}(\Ov{\vt}, \Grad \vU) - r \Grad G,
\]
we get
\begin{align} 
	&\intO{ \vre \left[ \frac{1}{\ep} \Grad G  - \partial_t \vU - (\vU \cdot \Grad) \vU \right] \cdot (\vue - \vU) } \br 
	&\quad 
	= 	\intO{ \frac{\vre}{\Ov{\vr}} \left[ \frac{1}{\ep} \Ov{\vr} \Grad G  + \Grad \Pi  - \Div \mathbb{S}(\Ov{\vt}, \Grad \vU) - r \Grad G \right] \cdot (\vue - \vU) }. 
	\nonumber
	\end{align}

Thus we can use the convergence established in \eqref{wc1}--\eqref{wc3} to rewrite \eqref{L6} in the form
 \begin{align}
	&\intO{ E_\ep \left(\vre, \vte, \vue \Big| \Ov{\vr} + \ep r, \Ov{\vt} + \ep \MTC, \vU \right)(\tau, \cdot) } \br 
	&+ \int_0^\tau \intO{ \Big( \mathbb{S} (\Ov{\vt}, \Ds \vue) - \mathbb{S} (\Ov{\vt}, \Ds \vU) \Big) : \Big( \Ds \vue - \Ds \vU  \Big) } \dt \br &
	+\int_0^\tau \intO{  \left(  \frac{\Ov{\vt} + \ep \MTC}{\vte^2} \right) \frac{\kappa (\vte) \Grad \vte \cdot \Grad \vte }{\ep^2}  } \dt
	- \int_0^\tau \intO{ \frac{\kappa (\vte) }{\vte}  \frac{\Grad \vte }{\ep} \cdot \Grad \MTC      } \dt \br 
	&\leq - \frac{1}{\ep} \int_0^\tau \intO{ \vre \Big[ (s(\vre, \vte) - s( \Ov{\vr} + \ep r,\Ov{\vt} + \ep \MTC ) \Big] \partial_t \MTC } \dt \br 
	&\quad - \frac{1}{\ep} \int_0^\tau \intO{ \vre \Big[ s(\vre, \vte) - s( \Ov{\vr} + \ep r, \Ov{\vt} + \ep \MTC) \Big] \vue \cdot \Grad \MTC  } \dt \br
	&\quad + \int_0^\tau \intO{ \frac{\vre}{\Ov{\vr}} \Grad \Pi \cdot (\vue - \vU) } \dt \br 
	&\quad + \int_0^\tau \intO{ \frac{\vre}{\Ov{\vr}} \left[ \frac{1}{\ep} \Ov{\vr} \Grad G    - r \Grad G \right] \cdot (\vue - \vU) } \dt \br 
	&\quad + \frac{1}{\ep^2} \int_0^\tau \intO{ \left[ \left( 1 - \frac{\vre}{\Ov{\vr} + \ep r} \right) \partial_t p(\Ov{\vr} + \ep r, \Ov{\vt} + \ep \MTC) - \frac{\vre}{\Ov{\vr} + \ep r} \vue \cdot \Grad p(\Ov{\vr} + \ep r, \Ov{\vt} + \ep \MTC) \right] } \dt\br 
	&\quad + C \int_0^\tau \intO{E_\ep \left(\vre, \vte, \vue \Big| \Ov{\vr} + \ep r, \Ov{\vt} + \ep \MTC, \vU \right) } \dt   + \mathcal{O}(\ep),
	\nonumber
\end{align}
where the symbol $\mathcal{O}(\ep)$ denotes a generic error, $\mathcal{O}(\ep) \to 0$ as $\ep \to 0$.
\textcolor{black}{Note that the convective term}
\[
\left| \intO{ \vre (\vue - \vc{U}) \otimes (\vue - \vc{U}) } \right| \leq 
\intO{ \vre |\vue - \vc{U} |^2 }
\]
\textcolor{black}{is controlled by the relative energy.}

In addition, in view of the convergences \eqref{wc1}--\eqref{wc3}, we conclude
\[
\int_0^\tau \intO{ \frac{\vre}{\Ov{\vr}} \Grad \Pi \cdot (\vue - \vU) } \dt = \mathcal{O}(\ep); 
\]
hence 
 \begin{align}
& \intO{ E_\ep \left(\vre, \vte, \vue \Big| \Ov{\vr} + \ep r, \Ov{\vt} + \ep \MTC, \vU \right) (\tau, \cdot) } \br 
&+ \int_0^\tau \intO{ \Big( \mathbb{S} (\Ov{\vt}, \Ds \vue) - \mathbb{S} (\Ov{\vt}, \Ds \vU) \Big) : \Big( \Ds \vue - \Ds \vU  \Big) } \dt \br &
+\int_0^\tau \intO{  \left(  \frac{\Ov{\vt} + \ep \MTC}{\vte^2} \right) \frac{\kappa (\vte) \Grad \vte \cdot \Grad \vte }{\ep^2}  } \dt
- \int_0^\tau \intO{ \frac{\kappa (\vte) }{\vte}  \frac{\Grad \vte }{\ep} \cdot \Grad \MTC      } \dt \br 
&\leq - \frac{1}{\ep} \int_0^\tau \intO{ \vre \Big[ (s(\vre, \vte) - s( \Ov{\vr} + \ep r,\Ov{\vt} + \ep \MTC ) \Big] \partial_t \MTC } \dt \br 
&\quad - \frac{1}{\ep} \int_0^\tau \intO{ \vre \Big[ s(\vre, \vte) - s( \Ov{\vr} + \ep r, \Ov{\vt} + \ep \MTC) \Big] \vue \cdot \Grad \MTC  } \dt \br
&\quad + \int_0^\tau \intO{ \frac{\vre}{\Ov{\vr}} \left[ \frac{1}{\ep} \Ov{\vr} \Grad G    - r \Grad G \right] \cdot (\vue - \vU) } \dt \br 
	&\quad + \frac{1}{\ep^2} \int_0^\tau \intO{ \left[ \left( 1 - \frac{\vre}{\Ov{\vr} + \ep r} \right) \partial_t p(\Ov{\vr} + \ep r, \Ov{\vt} + \ep \MTC) - \frac{\vre}{\Ov{\vr} + \ep r} \vue \cdot \Grad p(\Ov{\vr} + \ep r, \Ov{\vt} + \ep \MTC) \right] } \dt\br 
		&\quad + C \int_0^\tau \intO{E_\ep \left(\vre, \vte, \vue \Big| \Ov{\vr} + \ep r, \Ov{\vt} + \ep \MTC, \vU \right) } \dt   + \mathcal{O}(\ep).
	\label{L7}
\end{align}

\bigskip 

\noindent
{\bf Step 3:} At this stage, we use  the Boussinesq relation \eqref{OB1}
to obtain
\[
\frac{\partial p(\Ov{\vr}, \Ov{\vt} ) }{\partial \vr}  r + 
\frac{\partial p(\Ov{\vr}, \Ov{\vt} ) }{\partial \vt}  \MTC = \Ov{\vr}  G + \chi(t), 
\]
where 
\[
\chi = \frac{\partial p(\Ov{\vr}, \Ov{\vt} ) }{\partial \vt} \avintO{ \MTC } - \Ov{\vr} \avintO{ G } = \frac{\partial p(\Ov{\vr}, \Ov{\vt} ) }{\partial \vt} \avintO{ \MTC }
\]
{\cred since $\intO{ G }=0$.} Consequently, 
\begin{align}
\frac{1}{\ep^2} & \intO{ \left( 1 - \frac{\vre}{\Ov{\vr} + \ep r} \right) \partial_t p (\Ov{\vr} + \ep r, \Ov{\vt} + \ep \MTC) } \br &= 
\frac{1}{\ep} \intO{ \left( 1 - \frac{\vre}{\Ov{\vr} + \ep r} \right) \left( \frac{\partial p (\Ov{\vr} + \ep r, \Ov{\vt} + \ep \MTC) }{\partial \vr} \partial_t r +  \frac{\partial p (\Ov{\vr} + \ep r, \Ov{\vt} + \ep \MTC) }{\partial \vt} \partial_t \MTC  \right)	} \br 
&= \intO{ \frac{1}{\ep} \left( 1 - \frac{\vre}{\Ov{\vr} + \ep r} \right)  \left( \frac{\partial p(\Ov{\vr} + \ep r, \Ov{\vt} + \ep \MTC)}{\partial \vr} -\frac{\partial p (\Ov{\vr} , \Ov{\vt} ) }{\partial \vr} \right)  \partial_t r 
} \br &\quad + \intO{ \frac{1}{\ep}  \left( 1 - \frac{\vre}{\Ov{\vr} + \ep r} \right) 
	\left( \frac{\partial p (\Ov{\vr} + \ep r, \Ov{\vt} + \ep \MTC) }{\partial \vt} -\frac{\partial p (\Ov{\vr} , \Ov{\vt} ) }{\partial \vt} \right) \partial_t \MTC  	} \br
&\quad + { \frac{1}{\ep} \intO{ \left( \frac{\Ov{\vr} + \ep r - \vre}{\Ov{\vr} + \ep r} \right) \partial_t \chi }},
\label{L10}
	\end{align}
{where}
\begin{align}
&\intO{ \frac{1}{\ep} \left( 1 - \frac{\vre}{\Ov{\vr} + \ep r} \right)  \left( \frac{\partial p(\Ov{\vr} + \ep r, \Ov{\vt} + \ep \MTC)}{\partial \vr} -\frac{\partial p (\Ov{\vr} , \Ov{\vt} ) }{\partial \vr} \right)  \partial_t r 
} \br &\quad + \intO{ \frac{1}{\ep}  \left( 1 - \frac{\vre}{\Ov{\vr} + \ep r} \right) 
	\left( \frac{\partial p (\Ov{\vr} + \ep r, \Ov{\vt} + \ep \MTC) }{\partial \vt} -\frac{\partial p (\Ov{\vr} , \Ov{\vt} ) }{\partial \vt} \right) \partial_t \MTC  	} 
 = \mathcal{O}(\ep).
 \nonumber
\end{align}
{Moreover, }
\[
{ \frac{1}{\ep} \frac{\Ov{\vr} + \ep r - \vre}{\Ov{\vr} + \ep r} = -  \frac{\vre - \Ov{\vr}}{\ep( \Ov{\vr} + \ep r )} + 
\frac{r}{\Ov{\vr} + \ep r} \to \frac{1}{\Ov{\vr}} (r - \mathfrak{R} ).} 
\]
Seeing that
\[
\intO{ r } = \intO{ \mathfrak{R} } = 0,
\]
we may infer
\[
 \frac{1}{\ep} \intO{ \left( \frac{ \Ov{\vr} + \ep r - \vre}{\Ov{\vr} + \ep r } \right) \partial_t \chi } = \mathcal{O}(\ep).
\]

Similarly, 
\begin{align}
	&- \frac{1}{\ep^2}  \int_0^\tau \intO{ \frac{\vre}{\Ov{\vr} + \ep r} \vue \cdot \Grad p(\Ov{\vr} + \ep r, \Ov{\vt} + \ep \MTC) } \dt \br &= 
	- \frac{1}{\ep} \int_0^\tau \intO{ \frac{\vre}{\Ov{\vr} + \ep r} \vue \cdot \left( \frac{\partial p (\Ov{\vr} + \ep r, \Ov{\vt} + \ep \MTC) }{\partial \vr} \Grad r +  \frac{\partial p (\Ov{\vr} + \ep r, \Ov{\vt} + \ep \MTC) }{\partial \vt} \Grad \MTC  \right)	} \dt \br 
		&=- \int_0^\tau \intO{ \frac{1}{\ep} \frac{\vre}{\Ov{\vr} + \ep r} \vue \cdot \Grad r \left( \frac{\partial p (\Ov{\vr} + \ep r, \Ov{\vt} + \ep \MTC) }{\partial \vr} -\frac{\partial p (\Ov{\vr} , \Ov{\vt} ) }{\partial \vr} \right) } \dt \br
	&\quad - \int_0^\tau \intO{ \frac{1}{\ep} \frac{\vre}{\Ov{\vr} + \ep r} \vue \cdot \Grad \MTC \left( \frac{\partial p (\Ov{\vr} + \ep r, \Ov{\vt} + \ep \MTC) }{\partial \vt} -\frac{\partial p (\Ov{\vr} , \Ov{\vt} ) }{\partial \vt} \right) } \dt \br
	&- \frac{1}{\ep} \int_0^\tau \intO{\frac{\vre}{\Ov{\vr} + \ep r} \Ov{\vr} \vue \cdot \Grad G  } \dt .
	\label{L11}
\end{align}
Using \eqref{wc1}, \eqref{wc3}, we perform the limit in the first integral obtaining 
\begin{align} 
\int_0^\tau &\intO{ \frac{1}{\ep} \frac{\vre}{\Ov{\vr} + \ep r} \vue \cdot \Grad r \left( \frac{\partial p(\Ov{\vr} + \ep r, \Ov{\vt} + \ep \MTC) }{\partial \vr} -\frac{\partial p (\Ov{\vr} , \Ov{\vt} ) }{\partial \vr} \right) } \dt \br
& + \int_0^\tau \intO{ \frac{1}{\ep} \frac{\vre}{\Ov{\vr} + \ep r} \vue \cdot \Grad \MTC \left( \frac{\partial p (\Ov{\vr} + \ep r, \Ov{\vt} + \ep \MTC) }{\partial \vt} -\frac{\partial p (\Ov{\vr} , \Ov{\vt} ) }{\partial \vt} \right) } \dt \br
&= \int_0^\tau \intO{ \vu \cdot \left(  \frac{\partial^2 p (\Ov{\vr} , \Ov{\vt}) }{\partial^2 \vr} r \Grad r + \frac{\partial^2 p (\Ov{\vr} , \Ov{\vt}) }{\partial \vr \partial \vt }  \Grad (r \MTC) + 
\frac{\partial^2 p (\Ov{\vr} , \Ov{\vt}) }{\partial^2 \vt} {\cred \MTC} \Grad \MTC		\right) } + \mathcal{O}(\ep) \br 
&= \mathcal{O}(\ep)
	\nonumber
	\end{align}
as $\Div \vu = 0$.

Finally, 
\begin{align}
- \frac{1}{\ep} &\int_0^\tau \intO{\frac{\vre}{\Ov{\vr} + \ep r} \Ov{\vr} \vue \cdot \Grad G  } \dt\br & = 
- \frac{1}{\ep} \int_0^\tau \intO{\vre \vue \cdot \Grad G  } \dt + \frac{1}{\ep} \int_0^\tau \intO{ \left( 1 - \frac{\Ov{\vr}}{\Ov{\vr} + \ep r} \right) \vre \vue \cdot \Grad G  } \dt \br 
&= - \frac{1}{\ep} \int_0^\tau \intO{\vre \vue \cdot \Grad G  } \dt + \int_0^\tau \intO{ {\cred \frac{r}{\Ov{\vr} + \ep r}} \vre \vue \cdot \Grad G } \dt \br
&= - \frac{1}{\ep} \int_0^\tau \intO{\vre \vue \cdot \Grad G  } \dt + \int_0^\tau \intO{ r \vu \cdot \Grad G } \dt + \mathcal{O}(\ep).
\nonumber	
	\end{align}

Consequently, {\cred abbreviating $s_\ep=s(\vre, \vte)$,} relation \eqref{L7} can be rewritten as
 \begin{align}
&\intO{ E_\ep \left(\vre, \vte, \vue \Big| \Ov{\vr} + \ep r, \Ov{\vt} + \ep \MTC, \vU \right) (\tau, \cdot) }  \br 
&+ \int_0^\tau \intO{ \Big( \mathbb{S} (\Ov{\vt}, \Ds \vue) - \mathbb{S} (\Ov{\vt}, \Ds \vU) \Big) : \Big( \Ds \vue - \Ds \vU  \Big) } \dt \br &
+\int_0^\tau \intO{  \left(  \frac{\Ov{\vt} + \ep \MTC}{\vte^2} \right) \frac{\kappa (\vte) \Grad \vte \cdot \Grad \vte }{\ep^2}  } \dt
- \int_0^\tau \intO{ \frac{\kappa (\vte) }{\vte}  \frac{\Grad \vte }{\ep} \cdot \Grad \MTC      } \dt \br	
	&\leq - \frac{1}{\ep} \int_0^\tau \intO{ \left( \vre (s_\ep - s(\Ov{\vr} + \ep r, \Ov{\vt} + \ep \MTC)) \partial_t \MTC + \vre (s_\ep - s(\Ov{\vr} + \ep r, \Ov{\vt} + \ep \MTC)) \vue \cdot \Grad \MTC  \right) } \dt \br  
	&- \int_0^\tau \intO{ \frac{\vre}{\ep}  \Grad G \cdot \vU } \dt    +  \int_0^\tau \intO{ r \Grad G  \cdot \vU } \dt  \br 
&+ C \int_0^\tau \intO{E_\ep \left(\vre, \vte, \vue \Big| \Ov{\vr} + \ep r, \Ov{\vt} + \ep \MTC, \vU \right) } \dt   + \mathcal{O}(\ep).
	\label{L12}
\end{align}

\bigskip

{\bf Step 4:} In view of solenoidality $\Div \vU = 0$, we have
\begin{equation} \label{S11}
- \int_0^\tau \intO{ \frac{\vre}{\ep}  \Grad G \cdot \vU } \dt    +  \int_0^\tau \intO{ r \Grad G  \cdot \vU } \dt = 
- \int_0^\tau \intO{ \frac{\vre - (\Ov{\vr} + \ep r)}{\ep} \vU \cdot \Grad G } \dt.
\end{equation}
Consequently, we may use the bounds \eqref{be6}, \eqref{be7} along with the convergences established in \eqref{wc1}--\eqref{wc6}
to rewrite \eqref{L12} in the form 
 \begin{align}
&\intO{ E_\ep \left(\vre, \vte, \vue \Big| \Ov{\vr} + \ep r, \Ov{\vt} + \ep \MTC, \vU \right) (\tau, \cdot)}  \br 
&+ \int_0^\tau \intO{ \Big( \mathbb{S} (\Ov{\vt}, \Ds \vue) - \mathbb{S} (\Ov{\vt}, \Ds \vU) \Big) : \Big( \Ds \vue - \Ds \vU  \Big) } \dt \br &
+\int_0^\tau \intO{  \left(  \frac{\Ov{\vt} + \ep \MTC}{\vte^2} \right) \frac{\kappa (\vte) \Grad \vte \cdot \Grad \vte }{\ep^2}  } \dt
- \int_0^\tau \intO{ \frac{\kappa (\Ov{\vt}) }{\Ov{\vt}} \Grad \mathfrak{T} \cdot \Grad {\cred \MTC} } \dt \br	
&\leq - \frac{1}{\ep} \int_0^\tau \intO{ \vre \Big[ (s(\vre, \vte) - s( \Ov{\vr} + \ep r,\Ov{\vt} + \ep \MTC ) \Big] \partial_t \MTC } \dt \br 
&\quad - \frac{1}{\ep} \int_0^\tau \intO{ \vre \Big[ s(\vre, \vte) - s( \Ov{\vr} + \ep r, \Ov{\vt} + \ep \MTC) \Big] \vue \cdot \Grad \MTC  } \dt \br
&\quad+ \frac{1}{\ep} \int_0^\tau \intO{  \vre (s(\vre, \vte) - s(\Ov{\vr} + \ep r, \Ov{\vt} + \ep \MTC)) ( \vU - \vue) \cdot \Grad \MTC } \dt \br
	&\quad    +  \int_0^\tau \intO{ ( r - \mathfrak{R} ) \Grad G  \cdot \vU } \dt  \br 
&\quad+ C \int_0^\tau \intO{E_\ep \left(\vre, \vte, \vue \Big| \Ov{\vr} + \ep r, \Ov{\vt} + \ep \MTC, \vU \right) } \dt   + \mathcal{O}(\ep).
	\label{S16}
\end{align}

\bigskip 

{\bf Step 5:} Now we use the fact that $\MTC$ solves the modified heat equation \eqref{OBs}, specifically, 
\begin{align} 
	\partial_t \MTC + \vc{U} \cdot \Grad \MTC &= \frac{\Ov{\vt} \alpha (\Ov{\vr}, \Ov{\vt} ) }{ c_p (\Ov{\vr}, \Ov{\vt} )} \Grad G \cdot \vU + 
	\frac{\kappa(\Ov{\vt})}{\Ov{\vr} c_p (\Ov{\vr}, \Ov{\vt} )}  \Del \MTC 
		+ \frac{1}{\Ov{\vr} c_p (\Ov{\vr}, \Ov{\vt} )} \Lambda(t), \br \Lambda &= 
	\Ov{\vt} \alpha (\Ov{\vr}, \Ov{\vt}) \frac{\partial p (\Ov{\vr}, \Ov{\vt})}{\partial \vt} 
	\partial_t \avintO{ \MTC }.	\label{S17}
	\end{align}
Thus we may perform the limit in several integrals in \eqref{S16} obtaining
 \begin{align}
& \intO{ E_\ep \left(\vre, \vte, \vue \Big| \Ov{\vr} + \ep r, \Ov{\vt} + \ep \MTC, \vU \right) (\tau, \cdot) }  \br 
&+ \int_0^\tau \intO{ \Big( \mathbb{S} (\Ov{\vt}, \Ds \vue) - \mathbb{S} (\Ov{\vt}, \Ds \vU) \Big) : \Big( \Ds \vue - \Ds \vU  \Big) } \dt \br &
+\int_0^\tau \intO{  \left(  \frac{\Ov{\vt} + \ep \MTC}{\vte^2} \right) \frac{\kappa (\vte) \Grad \vte \cdot \Grad \vte }{\ep^2}  } \dt
- \int_0^\tau \intO{ \frac{\kappa (\Ov{\vt}) }{\Ov{\vt}} \Grad \mathfrak{T} \cdot \Grad \MTC      } \dt \br	
	&\leq -  \int_0^\tau \intO{  \Ov{\vr} \left(\frac{\partial s(\Ov{\vr}, \Ov{\vt})}{\partial \vr}(\mathfrak{R} - r) + \frac{\partial s(\Ov{\vr}, \Ov{\vt})}{\partial \vt}(\mathfrak{T} - \MTC) \right) {\frac{\Ov{\vt} \alpha (\Ov{\vr}, \Ov{\vt} )}{ c_p (\Ov{\vr}, \Ov{\vt} )} } \Grad G \cdot \vU                       } \dt \br
&\quad - { \int_0^\tau \intO{  \Ov{\vr} \left(\frac{\partial s(\Ov{\vr}, \Ov{\vt})}{\partial \vr}(\mathfrak{R} - r) + \frac{\partial s(\Ov{\vr}, \Ov{\vt})}{\partial \vt}(\mathfrak{T} - \MTC) \right)  \frac{\kappa(\Ov{\vt})}{\Ov{\vr} c_p (\Ov{\vr}, \Ov{\vt} )} \Del \MTC                  } \dt} \br	
&\quad { -  \int_0^\tau \intO{  \Ov{\vr} \left(\frac{\partial s(\Ov{\vr}, \Ov{\vt})}{\partial \vr}(\mathfrak{R} - r) + \frac{\partial s(\Ov{\vr}, \Ov{\vt})}{\partial \vt}(\mathfrak{T} - \MTC) \right)  \frac{1}{\Ov{\vr} c_p (\Ov{\vr}, \Ov{\vt} )} \Lambda(t)                  } \dt} \br
	& \quad   +  \int_0^\tau \intO{ ( r - \mathfrak{R} ) \Grad G  \cdot {\vc{U}} } \dt  \br 
&\quad + C \int_0^\tau \intO{E_\ep \left(\vre, \vte, \vue \Big| \Ov{\vr} + \ep r, \Ov{\vt} + \ep \MTC, \vU \right) } \dt   + \mathcal{O}(\ep),
	\label{S18}
\end{align}
where we have used 
\begin{align}
\frac{1}{\ep} &\int_0^\tau \intO{  \vre (s(\vr_\ep,\vt_\ep) - s(\Ov{\vr} + \ep r, \Ov{\vt} + \ep \MTC)) ( \vU - \vue) \cdot \Grad {\cred \MTC} } \dt \br
&\aleq \int_0^\tau \intO{E_\ep \left(\vre, \vte, \vue \Big| \Ov{\vr} + \ep r, \Ov{\vt} + \ep \MTC, \vU \right) } \dt.
\nonumber
\end{align}

Now, we use
\[
\intO{ (r - \mathfrak{R}) } = 0
\]
to rewrite the third integral on the right-hand side of \eqref{S18} as
\begin{align}
-  \int_0^\tau \intO{  \Ov{\vr} \left(\frac{\partial s(\Ov{\vr}, \Ov{\vt})}{\partial \vr}(\mathfrak{R} - r) + \frac{\partial s(\Ov{\vr}, \Ov{\vt})}{\partial \vt}(\mathfrak{T} - \MTC) \right)  \frac{1}{\Ov{\vr} c_p (\Ov{\vr}, \Ov{\vt} )} \Lambda(t)                  } \dt \nonumber \br 
=- \int_0^\tau \intO{  \frac{\partial s(\Ov{\vr}, \Ov{\vt})}{\partial \vt} \left[
	\frac{\partial p(\Ov{\vr}, \Ov{\vt} )}{\partial \vr} \left( \frac{\partial p(\Ov{\vr}, \Ov{\vt} )}{\partial \vt} \right)^{-1}(\mathfrak{R} - r) + (\mathfrak{T} - \MTC) \right]  \frac{1}{ c_p (\Ov{\vr}, \Ov{\vt} )} \Lambda(t)                  } \dt.
\label{S22}
\end{align}
In view of the Boussinesq relations \eqref{wc7}, \eqref{OB1}, the expression 
\[
\left[
\frac{\partial p(\Ov{\vr}, \Ov{\vt} )}{\partial \vr} \left( \frac{\partial p(\Ov{\vr}, \Ov{\vt} )}{\partial \vt} \right)^{-1}(\mathfrak{R} - r) + (\mathfrak{T} - \MTC) \right]
\]
is spatially homogeneous, meaning it depends on $t$ only.

Similarly, we can rewrite the second integral on the right-hand side of \eqref{S18} as 
\begin{align}
-\int_0^\tau \intO{  \Ov{\vr} \left(\frac{\partial s(\Ov{\vr}, \Ov{\vt})}{\partial \vr}(\mathfrak{R} - r) + \frac{\partial s(\Ov{\vr}, \Ov{\vt})}{\partial \vt}(\mathfrak{T} - \MTC) \right)  \frac{\kappa(\Ov{\vt})}{\Ov{\vr} c_p (\Ov{\vr}, \Ov{\vt} )} \Del \MTC                  } \dt \br 
= -	\int_0^\tau \intO{  \frac{\partial s(\Ov{\vr}, \Ov{\vt})}{\partial \vr}\left[ (\mathfrak{R} - r) +  \frac{\partial p(\Ov{\vr}, \Ov{\vt}) }{\partial \vt}
	\left(  \frac{\partial p(\Ov{\vr}, \Ov{\vt}) }{\partial \vr} \right)^{-1} (\mathfrak{T} - \MTC)  \right]  \frac{\kappa(\Ov{\vt})}{ c_p (\Ov{\vr}, \Ov{\vt} )} \Del \MTC                  } \dt \br 
+ \int_0^\tau \intO{ \left(  \frac{\partial s(\Ov{\vr}, \Ov{\vt})}{\partial \vr}\frac{\partial p(\Ov{\vr}, \Ov{\vt}) }{\partial \vt}
	\left(  \frac{\partial p(\Ov{\vr}, \Ov{\vt}) }{\partial \vr} \right)^{-1} (\mathfrak{T} - \MTC) - \frac{\partial s(\Ov{\vr}, \Ov{\vt})}{\partial \vt} (\mathfrak{T} - \MTC)\right) \frac{\kappa(\Ov{\vt})} { c_p (\Ov{\vr}, \Ov{\vt} )} \Del \MTC  } \dt.
\label{S23}	
	\end{align}
Similarly to the above, the quantity
\[
\left[ (\mathfrak{R} - r) +  \frac{\partial p(\Ov{\vr}, \Ov{\vt}) }{\partial \vt}
\left(  \frac{\partial p(\Ov{\vr}, \Ov{\vt}) }{\partial \vr} \right)^{-1} (\mathfrak{T} - \MTC)  \right]
\]
is independent of $x$.

Now, integrating equation \eqref{S17} in $x$  we obtain the identity
\begin{align} \left[ \left( \Ov{\vt} \alpha (\Ov{\vr} ,\Ov{\vt} ) 
	\frac{\partial p(\Ov{\vr}, \Ov{\vt})}{\partial \vt} \right)^{-1} -   \frac{1}{\Ov{\vr} c_p (\Ov{\vr}, \Ov{\vt}) } \right] |\Omega| \Lambda (t) = \int_{\partial \Omega} \frac{\kappa(\Ov{\vt})}{\Ov{\vr} c_p (\Ov{\vr}, \Ov{\vt} ) } \Grad \MTC \cdot \vc{n} \ \D \sigma_x.
	\label{S24}
	\end{align}

Finally, plugging \eqref{S24} into \eqref{S23} we can compute the sum of 
\eqref{S22} with the first integral in \eqref{S23} obtaining
\begin{align}
	-  &\intO{  \frac{\partial s(\Ov{\vr}, \Ov{\vt})}{\partial \vt} \left[
		\frac{\partial p(\Ov{\vr}, \Ov{\vt} )}{\partial \vr} \left( \frac{\partial p(\Ov{\vr}, \Ov{\vt} )}{\partial \vt} \right)^{-1}(\mathfrak{R} - r) + (\mathfrak{T} - \MTC) \right]  \frac{1}{ c_p (\Ov{\vr}, \Ov{\vt} )} \Lambda(t)                  } \br 
 -	 &\int_{\Omega}  \frac{\partial s(\Ov{\vr}, \Ov{\vt})}{\partial \vr}\left[ (\mathfrak{R} - r) +  \frac{\partial p(\Ov{\vr}, \Ov{\vt}) }{\partial \vt}
	\left(  \frac{\partial p(\Ov{\vr}, \Ov{\vt}) }{\partial \vr} \right)^{-1} (\mathfrak{T} - \MTC)  \right] \times \br &\times \left[ \Ov{\vr} \left( \Ov{\vt} \alpha (\Ov{\vr} ,\Ov{\vt} ) 
	\frac{\partial p(\Ov{\vr}, \Ov{\vt})}{\partial \vt} \right)^{-1} -   \frac{1}{ c_p (\Ov{\vr}, \Ov{\vt}) } \right] \Lambda (t) \dx 	.
	\label{S25}
	\end{align}

Now, in accordance with Gibbs' relation {\cred and the definitions of $\alpha$ and $c_p$ in \eqref{Coeff}}, 
\begin{align} 
	&- \frac{\partial s(\Ov{\vr}, \Ov{\vt})}{\partial \vt} 
\frac{\partial p(\Ov{\vr}, \Ov{\vt} )}{\partial \vr} \left( \frac{\partial p(\Ov{\vr}, \Ov{\vt} )}{\partial \vt} \right)^{-1} - \frac{\partial s(\Ov{\vr}, \Ov{\vt} )}{\partial \vr} 
\left( c_p (\Ov{\vr}, \Ov{\vt})\Ov{\vr} \left( \Ov{\vt} \alpha (\Ov{\vr} ,\Ov{\vt} ) 
\frac{\partial p(\Ov{\vr}, \Ov{\vt})}{\partial \vt} \right)^{-1} - 1 \right) \br 
= &- \frac{1}{\Ov{\vt}}\frac{\partial e(\Ov{\vr}, \Ov{\vt})}{\partial \vt} 
\frac{\partial p(\Ov{\vr}, \Ov{\vt} )}{\partial \vr} \left( \frac{\partial p(\Ov{\vr}, \Ov{\vt} )}{\partial \vt} \right)^{-1} \br &- \frac{\partial s(\Ov{\vr}, \Ov{\vt} )}{\partial \vr}
\left[ \Ov{\vr} \left( \frac{\partial e(\Ov{\vr}, \Ov{\vt}) }{\partial \vt} + \frac{1}{\Ov{\vr}} \Ov{\vt} \alpha (\Ov{\vr}, \Ov{\vt}) \frac{\partial p (\Ov{\vr}, \Ov{\vt} )}{\partial \vt} \right)\left( \Ov{\vt} \alpha (\Ov{\vr} ,\Ov{\vt} ) 
\frac{\partial p(\Ov{\vr}, \Ov{\vt})}{\partial \vt} \right)^{-1}
  - 1 \right] \br 
  = &- \frac{1}{\Ov{\vt}}\frac{\partial e(\Ov{\vr}, \Ov{\vt})}{\partial \vt} 
  \frac{\partial p(\Ov{\vr}, \Ov{\vt} )}{\partial \vr} \left( \frac{\partial p(\Ov{\vr}, \Ov{\vt} )}{\partial \vt} \right)^{-1}  
  - \frac{\partial s(\Ov{\vr}, \Ov{\vt} )}{\partial \vr}
  \left[ \Ov{\vr} \frac{\partial e(\Ov{\vr}, \Ov{\vt}) }{\partial \vt}\left(
  \Ov{\vt} \alpha (\Ov{\vr} ,\Ov{\vt} ) 
  \frac{\partial p(\Ov{\vr}, \Ov{\vt})}{\partial \vt} \right)^{-1} 
   \right] \br 
 = &- \frac{1}{\Ov{\vt}}\frac{\partial e(\Ov{\vr}, \Ov{\vt})}{\partial \vt} 
 \frac{\partial p(\Ov{\vr}, \Ov{\vt} )}{\partial \vr} \left( \frac{\partial p(\Ov{\vr}, \Ov{\vt} )}{\partial \vt} \right)^{-1}  
 + \frac{1}{\Ov{\vr}} \frac{\partial p(\Ov{\vr}, \Ov{\vt} )}{\partial \vt}
 \left[  \frac{\partial e(\Ov{\vr}, \Ov{\vt}) }{\partial \vt}\left(
 \Ov{\vt} \alpha (\Ov{\vr} ,\Ov{\vt} ) 
 \frac{\partial p(\Ov{\vr}, \Ov{\vt})}{\partial \vt} \right)^{-1} 
 \right] = 0.   
	\label{S26}
	\end{align}
Thus, the coefficient multiplying $\mathfrak{R} - r$ vanishes. By the same token, we deduce that the coefficient multiplying $\mathfrak{T} - \MTC$ vanishes.	

Next, we handle the second integral in \eqref{S23}. Using Gibbs' relation {\cred and the constitutive relations obtained in Section~\ref{EOS}}, specifically, 
\[
\frac{\partial s (\Ov{\vr}, \Ov{\vt} ) }{\partial \vr} = - \frac{1}{\Ov{\vr}^2} \frac{\partial p(\Ov{\vr}, \Ov{\vt})} 
{\partial \vt}.
\]
Thus, we get
\begin{align}  
\left[ \frac{\partial s(\Ov{\vr}, \Ov{\vt})}{\partial \vr}\frac{\partial p(\Ov{\vr}, \Ov{\vt}) }{\partial \vt}
	\left(  \frac{\partial p(\Ov{\vr}, \Ov{\vt}) }{\partial \vr} \right)^{-1}  - \frac{\partial s(\Ov{\vr}, \Ov{\vt})}{\partial \vt} \right]  \frac{\kappa(\Ov{\vt})} { c_p (\Ov{\vr}, \Ov{\vt} )} \br =
	-\left[ \frac{1}{\Ov{\vr}^2} \left( \frac{\partial p(\Ov{\vr}, \Ov{\vt}) }{\partial \vt} \right)^2
	\left(  \frac{\partial p(\Ov{\vr}, \Ov{\vt}) }{\partial \vr} \right)^{-1}  + \frac{1}{\Ov{\vt}} \frac{\partial e(\Ov{\vr}, \Ov{\vt})}{\partial \vt} \right]  \frac{\kappa(\Ov{\vt})} { c_p (\Ov{\vr}, \Ov{\vt} )} = {\cred -\frac{\kappa(\Ov{\vt})}{\Ov{\vt}}.}
\label{S27}
\end{align}

Finally, we regroup terms containing $\Grad G$: 
\begin{align}
-  \intO{  \Ov{\vr} \left(\frac{\partial s(\Ov{\vr}, \Ov{\vt})}{\partial \vr}(\mathfrak{R} - r) + \frac{\partial s(\Ov{\vr}, \Ov{\vt})}{\partial \vt}(\mathfrak{T} - \MTC) \right) {\frac{\Ov{\vt} \alpha (\Ov{\vr}, \Ov{\vt} )}{ c_p (\Ov{\vr}, \Ov{\vt} )} } \Grad G \cdot \vU                       }	%\br 
+   \intO{ ( r - \mathfrak{R} ) \Grad G  \cdot {\vc{U}} } \br 
=  \intO{  \Ov{\vr} \left(\frac{\partial s(\Ov{\vr}, \Ov{\vt})}{\partial \vr}\Grad (\mathfrak{R} - r) + \frac{\partial s(\Ov{\vr}, \Ov{\vt})}{\partial \vt}\Grad (\mathfrak{T} - \MTC) \right) {\frac{\Ov{\vt} \alpha (\Ov{\vr}, \Ov{\vt} )}{ c_p (\Ov{\vr}, \Ov{\vt} )} } G \cdot \vU                       } %\br
-   \intO{ \Grad ( r - \mathfrak{R} )  G  \cdot {\vc{U}} }.
	\label{S28}
	\end{align}

Using Boussinesq relation, we deduce
\begin{align}
&\Ov{\vr} \left(\frac{\partial s(\Ov{\vr}, \Ov{\vt})}{\partial \vr}\Grad (\mathfrak{R} - r) + \frac{\partial s(\Ov{\vr}, \Ov{\vt})}{\partial \vt}\Grad (\mathfrak{T} - \MTC) \right) {\frac{\Ov{\vt} \alpha (\Ov{\vr}, \Ov{\vt} )}{ c_p (\Ov{\vr}, \Ov{\vt} )} } \br 
=&\Ov{\vr} \left(\frac{\partial s(\Ov{\vr}, \Ov{\vt})}{\partial \vr}\Grad (\mathfrak{R} - r) - \frac{\partial s(\Ov{\vr}, \Ov{\vt})}{\partial \vt} {\cred \frac{\partial p(\Ov{\vr}, \Ov{\vt}) }{\partial \vr} 
\left( \frac{\partial p(\Ov{\vr}, \Ov{\vt}) }{\partial \vt} \right)^{-1} } \Grad (\mathfrak{R} - r) \right) {\frac{\Ov{\vt} \alpha (\Ov{\vr}, \Ov{\vt} )}{ c_p (\Ov{\vr}, \Ov{\vt} )} }	\br 
=& - \Ov{\vr} \left( \frac{1}{\Ov{\vr}^2} \frac{\partial p (\Ov{\vr}, \Ov{\vt}) }{\partial \vt} \Grad (\mathfrak{R} - r) + \frac{1}{\Ov{\vt}} \frac{\partial e(\Ov{\vr}, \Ov{\vt})}{\partial \vt} {\cred \frac{\partial p(\Ov{\vr}, \Ov{\vt}) }{\partial \vr} 
\left( \frac{\partial p(\Ov{\vr}, \Ov{\vt}) }{\partial \vt} \right)^{-1} }  \Grad (\mathfrak{R} - r) \right) {\frac{\Ov{\vt} \alpha (\Ov{\vr}, \Ov{\vt} )}{ c_p (\Ov{\vr}, \Ov{\vt} )} }\br 
=& - \Grad (\mathfrak{R} - r).
	\label{S29}
	\end{align}

Thus, {\cred rearranging terms and using $\mathfrak{T}|_{\partial\Omega}=\MTC|_{\partial\Omega}$,} \eqref{S18} reduces to the desired inequality 
\begin{align}
	&\intO{ E_\ep \left(\vre, \vte, \vue \Big| \Ov{\vr} + \ep r, \Ov{\vt} + \ep \MTC, \vU \right) (\tau, \cdot) }  \br 
	&+ \int_0^\tau \intO{ \Big( \mathbb{S} (\Ov{\vt}, \Ds \vue) - \mathbb{S} (\Ov{\vt}, \Ds \vU) \Big) : \Big( \Ds \vue - \Ds \vU  \Big) } \dt \br  &+ \int_0^\tau \intO{ \frac{\kappa (\Ov{\vt} ) }{\Ov{\vt}} \left|\Grad 
		\left(	\frac{\vte - \Ov{\vt}}{\ep} \right) - \Grad \MTC \right|^2 } \dt \br
	&\aleq \int_0^\tau \intO{E_\ep \left(\vre, \vte, \vue \Big|\Ov{\vr} + \ep r, \Ov{\vt} + \ep \MTC, \vU \right) } \dt   + \mathcal{O}(\ep).
	\label{S20}
\end{align}
Using Gr\" onwall's lemma and letting $\ep \to 0$ we obtain the conclusion claimed in Theorem \ref{MT1}.

%\bibliography{citace}

%\bibliographystyle{plain}

\def\cprime{$'$} \def\ocirc#1{\ifmmode\setbox0=\hbox{$#1$}\dimen0=\ht0
	\advance\dimen0 by1pt\rlap{\hbox to\wd0{\hss\raise\dimen0
			\hbox{\hskip.2em$\scriptscriptstyle\circ$}\hss}}#1\else {\accent"17 #1}\fi}

\end{document}